\documentclass[aos,MSNbibl,nameyear,seceqn,dvips]{arximspdf}
\usepackage{mathrsfs}
\usepackage{graphicx}
\doi{10.1214/13-AOS1086} 
\volume{41}
\issue{2}
\pubyear{2013}
\firstpage{568}
\lastpage{603}

\makeatletter
\newcommand{\rrVert}{\Vert}
\newcommand{\rrvert}{\vert}
\newcommand{\llVert}{\Vert}
\newcommand{\llvert}{\vert}
\newcommand{\eqref}[1]{(\ref{#1})}
\newcommand{\implies}{\Longrightarrow}
\newcommand{\mod}{\,\mathrm{mod}\,}
\newcommand{\iint}{\int\!\!\!\int}
\newtheorem{thmm}{Theorem}[section]
\newproclaim{rmk}[thmm]{Remark}
\newtheorem{cor}[thmm]{Corollary}
\newtheorem{prop}[thmm]{Proposition}
\newtheorem{lma}[thmm]{Lemma}
\newproclaim{condition}{Condition}
\renewcommand{\sc}[1]{{\langle #1 \rangle}}
\newcommand{\ud}{\,{d}}
\newcommand{\ee}{\mathbb{E}}
\newcommand{\eee}[1]{\mathbb{E}[#1]}

\newcommand{\cov}{\operatorname{cov}} 
\newcommand{\var}{\operatorname{var}} 

\newcommand{\vfi}{\varphi}
\newcommand{\vep}{\varepsilon}

\renewcommand{\b}[1]{{\mathbf{#1}}}
\newcommand{\bvep}{\bolds{\vep}}
\DeclareMathAlphabet{\mathpzc}{OT1}{pzc}{m}{it} 
\newcommand{\covK}{\mathpzc{r}} 
\newcommand{\covO}{\mathscr{R}} 
\newcommand{\lsub}[2]{{\vphantom{#2}}_{#1}{#2}} 

\newcommand{\udf}[2]{\widetilde{#1}_{#2}} 
\newcommand{\specdM}[1]{\Phi_{#1}} 
\newcommand{\specdK}[1]{\mathpzc{f}_{#1}} 
\newcommand{\specdKE}[1]{\hat{\mathpzc{f}}_{#1}} 
\newcommand{\specdO}[1]{\mathscr{F}_{#1}} 
\newcommand{\ii}{\mathrm{i}} 
\newcommand{\cum}[1]{\operatorname{cum}( #1 )} 
\newcommand{\hnorm}[1]{\Vert #1 \Vert} 
\newcommand{\snorm}[1]{\vert\!\vert\!\vert #1 \vert\!\vert\!\vert}

\makeatother

\begin{document}
\begin{frontmatter}

\title{Fourier analysis of stationary time series in function
space\thanksref{T1}}
\runtitle{Fourier analysis of functional time series}
\thankstext{T1}{Supported by a European Research Council Starting Grant Award.}

\begin{aug}
\author{\fnms{Victor M.} \snm{Panaretos}\corref{}\ead[label=e1]{victor.panaretos@epfl.ch}}
\and
\author{\fnms{Shahin} \snm{Tavakoli}\ead[label=e2]{shahin.tavakoli@epfl.ch}}
\runauthor{V.~M. Panaretos and S. Tavakoli}
\affiliation{Ecole Polytechnique F\'ed\'erale de Lausanne}
\address{Section de Math\'ematiques\\
Ecole Polytechnique F\'ed\'erale de Lausanne\\
1015 Lausanne\\
Switzerland\\
\printead{e1}\\
\phantom{E-mail:\ }\printead*{e2}}
\end{aug}

\received{\smonth{4} \syear{2012}}
\revised{\smonth{10} \syear{2012}}

%
\begin{abstract}
We develop the basic building blocks of a frequency domain framework
for drawing statistical inferences on the second-order structure of a
stationary sequence of functional data. The key element in such a
context is the spectral density operator, which generalises the notion
of a spectral density matrix to the functional setting, and
characterises the second-order dynamics of the process. Our main tool
is the functional Discrete Fourier Transform (fDFT). We derive an
asymptotic Gaussian representation of the fDFT, thus allowing the
transformation of the original collection of dependent random functions
into a collection of approximately independent complex-valued Gaussian
random functions. Our results are then employed in order to construct
estimators of the spectral density operator based on smoothed versions
of the periodogram kernel, the functional generalisation of the
periodogram matrix. The consistency and asymptotic law of these
estimators are studied in detail. As immediate consequences, we obtain
central limit theorems for the mean and the long-run covariance
operator of a stationary functional time series. Our results do not
depend on structural modelling assumptions, but only functional
versions of classical cumulant mixing conditions, and are shown to be
stable under discrete observation of the individual curves.
\end{abstract}

%


\end{frontmatter}

\section{Introduction}

In the usual context of functional data analysis, one wishes to make
inferences pertaining to the law of a continuous time stochastic
process $\{X(\tau);\tau\in[0,1]\}$ on the basis of a collection of $T$
realisations of this stochastic process, $\{X_t(\tau)\}_{t=0}^{T-1}$.
These are modelled as random elements of the separable Hilbert space
$L^2([0,1], \mathbb R)$ of square integrable real functions defined on
$[0,1]$. Statistical analyses typically focus on the first and
second-order characteristics of this law [see, e.g., \citet
{grenander1981}, \citet{rice1991}, \citet{ramsay2005}] and are, for
the most part, based on the fundamental Karhunen--Lo\`eve decomposition
[\citet{karhunen1947}, \citet{loeve1948}, \citet{dauxois1982}, \citet
{hall2006}]. Especially the second-order structure of random functions
is central to the analysis of functional data, as it is connected with
the smoothness properties of the random functions and their optimal
finite-dimensional representations [e.g., \citet{adler1990}]. When
functional data are independent and identically distributed, the entire
second-order structure is captured by the covariance operator [\citet
{grenander1981}], or related operators [e.g., \citet{locantore1999},
\citet{kraus2012}]. The assumption of identical distribution can be
relaxed, and this is often done by allowing a varying first-order
structure through the inclusion of covariate variables (or functions)
in the context of functional regression and analysis of variance
models; see \citet{cuevas2002}; \citet{cardot2003}; \citet{yao2005}.
Second-order structure has been studied in the ``nonidentically
distributed'' context mostly in terms of the so-called common principal
components model [e.g., \citet{benko2009}], in a comparison setting,
where two functional populations are compared with respect to their
covariance structure {[e.g., \citet{panaretos2010,
Boente2011,Horvath2012, Fremdt2011}], and in the context of
detection of sequential changes in the covariance structure of
functional observations [\citet{Horvath10a}]; see \citet{Horvath2012}
for an overview.

For sequences of potentially dependent functional data,
Gabrys and Ko\-koszka (\citeyear{Gabrys2007}) and \citet{Gabrys2010} study the detection of
correlation. To obtain a complete description of the second-order
structure of dependent functional sequences, one needs to consider
autocovariance operators relating different lags of the series, as is
the case in multivariate time series. This study will usually be
carried out under the assumption of stationarity. Research in this
context has mostly focused on stationary functional series that are
linear. Problems considered include that of the estimation of the
second-order structure [e.g., \citet{mas2000}, \citet{bosq2002},
\citet{dehling2005}] and that of prediction [e.g.,
\citet{antoniadis2003}, \citet{ferraty2004}, \citet{antoniadis2006}]. It
can be said that the linear case is now relatively well understood, and
\citet{Bosq00} and \citet{Blanke2008} provide a detailed overview
thereof.\looseness=-1

Recent work has attempted to move functional time series beyond linear
models and construct inferential procedures for time series that are
not a priori assumed to be described by a particular model, but are
only assumed to satisfy certain weak dependence conditions. \citet
{Hormann10} consider the effect that weak dependence can have on the
principal component analysis of functional\vadjust{\goodbreak} data and propose weak
dependence conditions under which they study the stability of
procedures that assume independence. They also study the problem of
inferring the long-run covariance operator by means of
finite-dimensional projections. \citet{Horvath11} give a central limit
theorem for the mean of a stationary weakly dependent functional
sequence, and propose a consistent estimator for the long-run
covariance operator.

In this paper, rather than focus on isolated characteristics such as
the long-run covariance, we consider the problem of inferring the
complete second-order structure of stationary functional time series
without any structural modelling assumptions, except for cumulant-type
mixing conditions. Our approach is to study the problem via Fourier
analysis, formulating a frequency domain framework for weakly dependent
functional data. To this aim, we employ suitable generalisations of
finite-dimensional notions [e.g., \citet{Brill01}, \citet
{bloomfield2000}, \citet{priestley2001}] and provide conditions for
these to be well defined.

We encode the complete second-order structure via the spectral density
operator, the Fourier transform of the collection of autocovariance
operators, seen as operator-valued functions of the lag argument; see
Proposition~\ref{lmaprop-specdk}. We propose strongly consistent and
asymptotically Gaussian estimators of the spectral density operator
based on smoothing the periodogram operator---the functional analogue
of the periodogram matrix; see Theorems~\ref{thmimse-convergence} and
\ref{thmCLTspecdensity}. In this sense, our methods can be seen as
functional smoothing, as overviewed in \citet{ferraty2006}, but in an
operator context; see also, for example, \citet{ferraty-goia2011},
\citet{ferraty-tadj2011}, \citet{louani2010}. As a by-product, we
also obtain central limit theorems for both the mean and long-run
covariance operator of stationary time series paralleling or extending
the results of \citet{Horvath11}, but under different weak dependence
conditions; see Corollaries~\ref{corCLT} and~\ref{corlimitlongrun}.
The key result employed in our analysis is the asymptotic
representation of the discrete Fourier transform of a weakly dependent
stationary functional process as a collection of independent Gaussian
elements of $ L^2( [0,1], \mathbb{C} )$, the Hilbert space of square integrable
complex-valued functions, with mean zero and covariance operator
proportional to the spectral density operator at the corresponding
frequency (Theorem~\ref{thmudf-asymptotics}). Weak dependence
conditions required to yield these results are moment type conditions
based on cumulant kernels, which are functional versions of cumulant
functions. A noteworthy feature of our results and methodology is that
they do not require the projection onto a finite-dimensional subspace,
as is often the case with functional time series [\citet{Hormann10},
\citet{Sen10}]. Rather, our asymptotic results hold for purely
infinite-dimensional functional data.

The paper is organised in seven sections and the supplementary material [\citet{supp}]. The
building blocks of the frequency domain framework are developed in
Section~\ref{secspectraldensity}. After some basic definitions and
introduction of notation, Section~\ref{subsecspec} provides
conditions for the definition of the spectral density operator. The
functional version of the discrete Fourier transform is introduced in
Section~\ref{secfDFT}, where its analytical and asymptotic
properties are investigated. Section~\ref{secperiodogram} then
introduces the periodogram operator and studies its mean and covariance
structure. The estimation of the spectral density operator by means of
smoothing is considered in Section~\ref{secSpecdensityestimation}.
Section~\ref{examples} provides a detailed discussion on the weak
dependence conditions introduced in earlier sections. The effect of
observing only discretely sampled functions is considered in Section
\ref{discrete}, where the consistency is seen to persist under
conditions on the nature of the discrete sampling scheme. Finite-sample
properties are illustrated via simulation in Section~\ref{simulations}.
Technical background and several lemmas required for the proofs or the
main results are provided in a an extensive supplementary material [\citet{supp}]. One of
our technical results, Lemma~\ref{lmatech-criterion-unif-tightness},
collects some results that may be of independent interest in functional
data analysis when seeking to establish tightness in order to extend
finite-dimensional convergence results to infinite dimensions; it is
given in the main paper, in a short section (Section~\ref{technical}).

\section{Spectral characteristics of stationary functional data}\label
{secspectraldensity}

We start out this section with an introduction of some basic
definitions and notation. Let $\{X_t\}_{t\in\mathbb{Z}}$ be a
functional time series indexed by the integers, interpreted as time.
That is, for each $t$, we understand $X_{t}$ as being a random element
of $ L^2( [0,1], \mathbb{R} )$, with
\[
\tau\mapsto X_{t}(\tau) \in\mathbb{R}, \qquad \tau\in[0,1],
\]
denoting its parametrisation. Though all our results will be valid for
any separable Hilbert space, we choose to concentrate on $ L^2( [0,1], \mathbb{R} )$, as
this is the paradigm for functional data analysis. We denote the inner
product in $ L^2( [0,1], \mathbb{R} )$ by $\langle\cdot,\cdot\rangle$, and the induced norm
by $\|\cdot\|_2$,
\[
\langle f,g\rangle= \biggl(\int_0^1f(\tau)g(
\tau)\ud\tau \biggr)^{1/2},\qquad \|g\|_2=\langle g,g
\rangle^{1/2},\qquad f,g\in L^2\bigl( [0,1], \mathbb{R}
\bigr).
\]
Equality of $L^{2}$ elements will be understood in the sense of the
norm of their difference being zero. The imaginary number will de
denoted by $\ii$, $\ii^2=-1$, and the complex conjugate of $z\in\mathbb
{C}$ will be denoted as $\bar{z}$. We also denote $\Delta^{(T)}(\omega)
= \sum_{t=0}^{T-1} \exp(-\ii\omega t).$ The Hermitian adjoint of an
operator $\mathscr{A}$ will be denoted as $\mathscr{A}^{\dag}$. For a
function $g\dvtx D \subset\mathbb R^n \rightarrow\mathbb C,$ we denote
$\hnorm{g}_\infty= \sup_{\b{x} \in D} |g(\b{x})|.$

Throughout, we assume that the series $\{X_t\}_{t\in\mathbb{Z}}$ is
strictly stationary: for any finite set of indices $I \subset\mathbb
Z$ and any $s \in\mathbb Z$, the joint law of $\{ X_{t}, t \in I\}$
coincides with that of $\{ X_{t + s}, t \in I\}$. If $\ee\hnorm
{X_{0}}_2 < \infty,$ the mean of $X_{t}$ is well defined, belongs to
$ L^2( [0,1], \mathbb{R} )$, and is independent of $t$ by stationarity,
$\mu(\tau) = \ee X_{t}(\tau).$
We also define the autocovariance kernel at lag $t$ by
\[
\covK_{t}(\tau, \sigma) = \mathbb{E} \bigl[ \bigl( X_{t+s}(
\tau) - \mu(\tau) \bigr) \bigl( X_{s}(\sigma) - \mu(\sigma) \bigr)
\bigr], \qquad\tau, \sigma\in [0,1] \mbox{ and } t,s \in\mathbb{Z}.
\]
This kernel is well defined in the $L^{2}$ sense if $\ee\hnorm
{X_{0}}^{2}_2 < \infty$; if continuity in mean square of $X_t$ is
assumed, then it is also well defined pointwise.
Each kernel $\covK_{t}$ induces a corresponding operator $\covO_t\dvtx  L^2( [0,1], \mathbb{R} )
\rightarrow L^2( [0,1], \mathbb{R} )$ by right integration, the \emph{autocovariance
operator at lag $t$},
\[
\covO_{t}h (\tau) = \int_0^1
\covK_{t} (\tau, \sigma) h (\sigma) \ud \sigma=\cov \bigl[\bigl
\langle X_0,h \rangle, X_t (\tau) \bigr], \qquad h \in L^2
\bigl( [0,1], \mathbb{R} \bigr).
\]
One of the notions we will employ to quantify the weak dependence among
the observations $\{X_{t}\}$ is that of a \emph{cumulant kernel} of the
series; the pointwise definition of a \emph{$k$\textup{th} order cumulant
kernel} is
\[
\operatorname{cum} \bigl( X_{t_{1}}(\tau_{1}), \ldots,
X_{t_{k}}(\tau_{k}) \bigr) = \sum_{\nu=
(\nu_{1}, \ldots, \nu_{p})}
(-1)^{p-1} (p-1)! \prod_{l=1}^{p}
\mathbb{E} \biggl[\prod_{j \in\nu_{l}} X_{t_{j}}(
\tau_{j}) \biggr],
\]
where the sum extends over all unordered partitions of $\{1, \ldots, k\}
$. 
Assuming
$\ee\hnorm{X_{0}}_{2}^{l} < \infty$ for $l\geq1$ guarantees that the
cumulant kernels are well defined {in an $L^2$ sense.}
A cumulant kernel of order $2k$ gives rise to a corresponding \emph
{$2k$\textup{th} order cumulant operator} $\mathscr{R}_{t_1,\ldots,t_{2k-1}}\dvtx L^2([0,1]^k,\mathbb{R})\rightarrow L^2([0,1]^k,\mathbb{R})$, defined by
right integration,
\begin{eqnarray*}
&&\mathscr{R}_{t_1,\ldots,t_{2k-1}} h(\tau_1,\ldots,\tau_k)\\
&&\qquad=\int
_{[0,1]^k} \operatorname{cum} \bigl(X_{t_{1}}(
\tau_{1}), \ldots, X_{t_{2k-1}}(\tau_{2k-1}),X_{0}(
\tau _{2k}) \bigr)\\
&&\hspace*{59pt}{}\times h( \tau_{k+1},\ldots,\tau_{2k})\ud\tau_{k+1}
\cdots\ud \tau_{2k}.
\end{eqnarray*}

\subsection{The spectral density operator}\label{subsecspec} The
autocovariance operators encode all the second-order dynamical
properties of the series and are typically the main focus of functional
time series analysis. Since we wish to formulate a framework for a
frequency domain analysis of the series $\{X_t\}$, we need to consider
a suitable notion of Fourier transform of these operators.
This we call the \emph{spectral density operator} of $\{X_t\}$, defined
rigorously in Proposition~\ref{lmaprop-specdk} below. {Results of a
similar flavour related to Fourier transforms between general Hilbert
spaces can be traced back to, for example, \citet{kolmogorov1978}; we
give here the precise versions that we will be requiring, for
completeness, since those results do not readily apply in our setting. }
%
\begin{prop} \label{lmaprop-specdk}
Suppose $p=2$ or $p=\infty,$ and consider the following conditions:

\textup{I}($p$) the autocovariance kernels satisfy $\sum_{t \in\mathbb Z}
\| \covK_{t} \|_{p} < \infty;$

\textup{II} the autocovariance operators satisfy $\sum_{t \in\mathbb Z}
\snorm{\covO_{t}}_{1} < \infty,$
where $\snorm{\covO_{t}}_{1} $ is the nuclear norm or Schatten 1-norm;\vadjust{\goodbreak}
see Paragraph~\textup{F.1.1} in the supplementary material [\citet{supp}]. Then,
under \textup{{I($p$)}}, for any $\omega\in\mathbb R$, the following
series converges in $\hnorm{\cdot}_p$:
%
\begin{equation}
\label{eqspecdk-lma} \specdK{\omega}(\cdot, \cdot) = \frac{1}{2\pi} \sum
_{t \in\mathbb Z} \exp(-\ii\omega t) \covK_{t}(\cdot,
\cdot).
\end{equation}
We call the limiting kernel $\specdK{\omega}$ the \emph{spectral
density kernel} at frequency $\omega$. It is uniformly bounded and also
uniformly continuous in $\omega$ with respect to $\hnorm{\cdot}_p$;
that is, given $\vep> 0$, there exists a $\delta> 0$ such that
\[
|\omega_1 - \omega_2| < \delta\quad\implies\quad\Vert \specdK{
\omega_1} - \specdK{\omega_2} \Vert_p < \vep.
\]
The spectral density operator $\specdO{\omega}$, the operator induced
by the spectral density kernel through right-integration, is
self-adjoint and nonnegative definite for all $\omega\in\mathbb R$.
Furthermore, the following inversion formula holds in $\hnorm{\cdot}_p$:
%
\begin{equation}
\label{eqinversion-fourier} \int_{0}^{2\pi}
\specdK{\alpha}(\tau, \sigma) e^{\ii t \alpha} \ud \alpha= \covK_{t}(
\tau, \sigma)\qquad\forall t, \tau, \sigma.
\end{equation}

Under only \textup{{II}}, we have
%
\begin{equation}
\label{eqlmaprop-specdkeq-0} \specdO{\omega} = \frac{1}{2\pi}\sum
_{t \in\mathbb Z} e^{-\ii\omega t} \covO_{t},
\end{equation}
where the convergence holds in nuclear norm. In particular, the
spectral density operators are nuclear, and
$\snorm{ \specdO{\omega} }_{1} \leq\frac{1}{2\pi}\sum_{t} \snorm{\covO_{t}}_{1} <
\infty.$
\end{prop}
\begin{pf}
See Proposition~A.1 in the supplementary material [\citet{supp}].
\end{pf}

The inversion relationship \eqref{eqinversion-fourier}, in particular,
shows that the autocovariance operators and the spectral density
operators comprise a Fourier pair, thus reducing the study of
second-order dynamics to that of the study of the spectral density operator.

{We use the term \emph{spectral density operator} by analogy to the
multivariate case, in which the Fourier transform of the
autocovariance functions is called the spectral density matrix; see,
for example, \citet{Brill01}. In our case, since the time series takes
values in
$ L^2( [0,1], \mathbb{R} )$, the autocovariance functions are in fact \emph{operators} and
their Fourier transform is an \emph{operator}, hence the term
\emph{spectral density operator}. In light of the inversion formula
\eqref{eqinversion-fourier}, for fixed
$(\tau, \sigma),$ we can think of the $\omega\mapsto\specdK{\omega
}(\tau, \sigma)$ as being a (complex) measure, giving the \emph
{distribution of energy between $X_t(\tau)$ and $X_0(\sigma)$ across
different frequencies.} That is, $\omega\mapsto\specdK{\omega}(\tau,
\tau)\geq0$ gives the power spectrum of the univariate time series $\{
X_t(\tau)\}_{t\in\mathbb{Z}}$, while given $\tau\neq\sigma$, $\omega
\mapsto\specdK{\omega}(\tau, \sigma)\in\mathbb{C}$ gives the
cross spectrum of the univariate time series $\{X_t(\tau)\}_{t\in\mathbb
{Z}}$ with $\{X_t(\sigma)\}_{t\in\mathbb{Z}}$. When a point-wise
interpretation of $\{X_t\}_{t\in\mathbb{Z}}$ is not possible (e.g.,
because it is only interpretable via $L^2$ equivalence classes), the
spectral density operator admits a weak interpretation as follows:
given $L^2$ elements $\psi\neq\phi$, the mapping $\omega\mapsto
\langle\psi,\specdO{\omega}\psi\rangle\geq0$ is the power spectrum of
the univariate time series $\{\langle\psi,X_t\rangle\}_{t\in\mathbb
{Z}}$, while $\omega\mapsto\langle\psi,\specdO{\omega}\phi\rangle=
\langle\specdO{\omega}\psi,\phi\rangle\in\mathbb{C}$ is the cross
spectrum of the univariate time series $\{\langle\psi,X_t\rangle\}
_{t\in\mathbb{Z}}$ with the univariate time series $\{\langle\phi,X_t\rangle\}_{t\in\mathbb{Z}}$. In this sense, $\specdO{\omega}$
provides a complete characterisation of the second-order dynamics of
the functional process $\{X_t\}$; {see also \citet{stoc} for the role
of the spectral density operator in the spectral representation and the
harmonic principal component analysis of functional time series}.}

\subsection{The functional discrete Fourier transform and its
properties}\label{secfDFT}

In practice, a stretch of length $T$ of the series $\{X_t\}_{t\in\mathbb
{Z}}$ will be available, and we will wish to draw inferences on the
spectral density operator based on this finite stretch. The main tool
that we will employ is the functional version of the discrete Fourier
transform (DFT). In particular, define the \emph{functional Discrete
Fourier Transform (fDFT)} of $\{X_t\}_{t=0}^{T-1}$ to be
\[
\udf{X} {\omega}^{(T)}(\tau) = (2\pi T)^{-1/2} \sum
_{t=0}^{T-1} X_{t}(\tau) \exp(-\ii\omega t).
\]
It is of interest to note here that the construction of the fDFT does
not require the representation of the data in a particular basis. The
fDFT transforms the $T$ functional observations to a mapping from
$\mathbb{R}$ into $ L^2( [0,1], \mathbb{C} )$. It straightforwardly inherits some basic
analytical properties that its finite-dimensional counterpart
satisfies; for example, it is $2\pi$-periodic and Hermitian with
respect to $\omega$, and linear with respect to the series $\{X_t\}$.

The extension of the stochastic properties of the multivariate DFT to
the fDFT, however, is not as straightforward. It is immediate that $\ee
\hnorm{\udf{X}{\omega}^{(T)}}_2^{l}<~\infty$ {if $\ee\hnorm{X_t}_2^l <
\infty$},
{and hence the fDFT is almost surely in $ L^2( [0,1], \mathbb{C} )$ if $\ee\hnorm{X_t}_2^2
< \infty$}. We will see that the asymptotic covariance operator of this
object coincides with the spectral density operator. Most importantly,
we prove below that the fundamental stochastic property of the
multivariate DFT can be adapted and extended to the
infinite-dimensional case; that is, under suitable weak dependence
conditions, as $T\rightarrow\infty$, the fDFT evaluated at distinct
frequencies yields independent and Gaussian random elements of $ L^2( [0,1], \mathbb{C} )$.
The important aspect of this limit theorem is that it does not require
the assumption of any particular model for the stationary series, and
imposes only cumulant mixing conditions. A more detailed discussion of
these conditions is provided in Section~\ref{examples}.

\begin{thmm}[(Asymptotic distribution of the fDFT)] \label{thmudf-asymptotics}
Let $\{X_t\}_{t=0}^{T-1}$ be a strictly stationary sequence of random
elements of $ L^2( [0,1], \mathbb{R} )$, of length~$T$. Assume the following conditions hold:
\begin{longlist}
\item[(i)] {$\ee\hnorm{X_0}_{2}^k < \infty$, }
$ \sum_{t_{1}, \ldots, t_{k-1} = -\infty}^{\infty} \llVert
\cum{X_{t_{1}}, \ldots, X_{t_{k-1}}, X_{0}} \rrVert _{2} < \infty,\break\forall k\geq2$.
\item[(ii)] $ \sum_{t \in\mathbb Z} \snorm{\covO_{t}}_{1}
< \infty.$
\end{longlist}

Then, for
$ \omega_{1,T}:= \omega_{1} = 0,$
$\omega_{2,T}:= \omega_{2} = \pi,$
and \emph{distinct} integers
\[
s_{3,T}, \ldots, s_{J,T} \in \bigl\{ 1, \ldots, \bigl
\lfloor(T-1)/2 \bigr\rfloor \bigr\}
\]
such that
\[
\omega_{j, T}:= \frac{2 \pi s_{j,T}}{T} \stackrel{T \rightarrow\infty } {
\longrightarrow} \omega_{j},\qquad j = 3, \ldots, J,
\]
we have
%
\begin{equation}
\label{eqdftclt} \udf{X} {\omega_{1}}^{(T)} - \sqrt{
\frac{T}{2\pi}} \mu \stackrel {d} {\longrightarrow} \udf{X} {
\omega_{1}}\qquad \mbox{as } T \rightarrow\infty,
\end{equation}
and
$\udf{X}{\omega_{j,T}}^{(T)} \stackrel{d}{\longrightarrow} \udf
{X}{\omega_{j}},  \mbox{as } T \rightarrow\infty,   j = 2, \ldots, J$
where $  \{\udf{X}{\omega_{j}} \}$ are independent mean zero
Gaussian elements of $L^2  ( [0,1], \mathbb{R}  )$ for
$j=1,2$, and of $L^2  ( [0,1], \mathbb{C}  )$ for $j = 3,
\ldots, J,$ with covariance operators $\specdO{\omega_{j}}$, respectively.
\end{thmm}

\begin{rmk}
Though the $\{\omega_{j,T}\}_{j=3}^{J}$ are distinct for every $T$,
the limiting frequencies $\{\omega_j\dvtx j=3,\ldots,J\}$ need not be distinct.
\end{rmk}

Note here that condition (i) with $k=2$ is already
required in order to define the spectral density kernel and operator in
Proposition~\ref{lmaprop-specdk}. Condition (i) for $k\geq3$ is
the generalisation of the standard multivariate cumulant condition to
the functional case
[\citet{Brill01}, Condition 2.6.1], and reduces to that exact same
condition if the data are finite-dimensional.
Condition (ii) is required so that the spectral
density operator be a nuclear operator at each $\omega$ [which is in
turn a necessary condition for the weak limit of the fDFT to be almost
surely in $ L^2( [0,1], \mathbb{C} )$]. As we shall see, condition (ii) is, in fact, a
sufficient condition for tightness of the fDFT, seen as a functional
process indexed by frequency.

\begin{pf*}{Proof of Theorem~\ref{thmudf-asymptotics}}
Consider $p_{\omega}^{(T)}(\tau, \sigma) = \udf{X}{\omega}^{(T)}(\tau)
\udf{X}{-\omega}^{(T)}(\sigma),$
and assume initially that $\mu= 0$. We will treat the case $\mu\neq
0$ at the end of the proof. First we show that for any $\omega$ (or
sequence $\omega_{T}$), the sequence of random elements $\udf{X}{\omega
}^{(T)}, T=1,2, \ldots,$ is tight. To do this, we shall use Lemma~\ref
{lmatech-criterion-unif-tightness}. Fix an orthonormal basis $\{\vfi
_{n}\}_{n\geq1}$ of $ L^2( [0,1], \mathbb{R} )$ and let $H= L^2( [0,1], \mathbb{C} )$. We notice that
$p^{(T)}_{\omega}$ is a random element of the (complete) tensor product
space $H \otimes H$, with scalar product and norm $\sc{\cdot, \cdot
}_{H\otimes H}, \hnorm{\cdot}_{H \otimes H}$, respectively; see \citeauthor{Weid80} [(\citeyear{Weid80}), Paragraph 3.4],
for instance. Notice that $\llvert \sc{\udf
{X}{\omega}, \vfi_n}\rrvert ^{2} = \sc{p_{\omega}^{(T)}, \vfi_n \otimes
\vfi_n}.$
Since $\ee\hnorm{p_{\omega}^{(T)}}_{H \otimes H} < \infty$ and the\vadjust{\goodbreak}
projection $\mathcal{P}_{n}\dvtx H \otimes H \rightarrow\mathbb C$ defined by
$\mathcal{P}_{n}(f) = \sc{f, \vfi_n \otimes\vfi_n}_{H \otimes H}$
is continuous and linear, we deduce
\begin{eqnarray*}
\ee \bigl\llvert
\bigl\langle \udf{X} {\omega}^{(T)}, \vfi_{n} \bigr\rangle
 \bigr
\rrvert ^{2} &=& \ee \mathcal{P}_{n}p_{\omega}^{(T)}
= \mathcal{P}_{n} \ee p_{\omega}^{(T)} =
\frac{1}{2\pi} \mathcal{P}_{n} \int_{-\pi}^{\pi}
F_{T}(\omega- \alpha)\specdK{\alpha} \,d\alpha
\\
&=& \frac{1}{2\pi} \int_{-\pi}^{\pi}
F_{T}(\omega- \alpha) \mathcal {P}_{n} \specdK{\alpha} \,d
\alpha\leq\sup_{\alpha\in\mathbb R} \llvert \mathcal{P}_{n}
\specdK{\alpha} \rrvert.
\end{eqnarray*}
The third equality comes from Proposition \ref
{propperiodogram-exact-mean} (which is independent of previous
results), the fourth equality follows from Tonelli's theorem [\citet{Wheeden1977}, page 92], and the last inequality is Young's inequality
[\citet{Hunter2005}, Theorem 12.58].
Notice that $ \llvert  \mathcal{P}_{n} \specdK{\alpha}\rrvert  = | \sc
{\specdO{\alpha} \vfi_{n}, \vfi_{n}}| \leq\sum_{t} | \sc{\covO_{t} \vfi
_{n}, \vfi_{n}} |$ by \eqref{eqlmaprop-specdkeq-0}.
Setting $a_{n} = \sum_{t} | \sc{\covO_{t} \vfi_{n}, \vfi_{n}} |$, which
is independent of $\alpha$ and $T$, we have
$\ee\llvert  \sc{\udf{X}{\omega}^{(T)}, \vfi_{n}} \rrvert ^{2} \leq a_{n},$
and $\sum_{n} a_{n} \leq\sum_{t \in\mathbb Z} \snorm{\covO_{t}}_{1} <
\infty.$
Therefore, we have proven that $\udf{X}{\omega}^{(T)}$ is tight.
Consequently, the random element
$ ( \udf{X}{\omega_{1,T}}^{(T)}, \ldots, \udf{X}{\omega
_{J,T}}^{(T)}  )$
of $ ( L^2  ( [0,1], \mathbb{C}  )  )^J$ is also tight.
Its asymptotic distribution is therefore determined by the convergence
of its finite-dimensional distributions; see, for example, \citeauthor{Ledo91} [(\citeyear{Ledo91}), Paragraph 2.1].
Thus, to complete the proof, it suffices to show that for any $\psi
_{1}, \ldots, \psi_{J} \in L^2( [0,1], \mathbb{C} )$,
%
\begin{equation}
\label{eqthmudf-asymptoticseq-0} \bigl(
\bigl\langle \udf{X} {\omega_{1,T}}^{(T)},
\psi_{1} \bigr\rangle, \ldots,
\bigl\langle \udf {X} {\omega_{J,T}}^{(T)},
\psi_{J} \bigr\rangle
 \bigr) \stackrel{d} {\longrightarrow } \bigl(
\langle \udf{X} {
\omega_{1}}, \psi_{1} \rangle, \ldots,
\langle \udf{X} {\omega
_{J}}, \psi_{J} \rangle
 \bigr),
\end{equation}
where $\udf{X}{\omega_{j}} \sim
{\mathcal{N}(0, \specdO{\omega_{j}} )}$ are
independent Gaussian random elements of $H$, where $H =  L^2( [0,1], \mathbb{R} )$ if $j=1,2$ and
$H= L^2( [0,1], \mathbb{C} )$ if $j = 3, \ldots, J$.
This is a consequence of the following claim, which is justified by
\citeauthor{Brill01} [(\citeyear{Brill01}), Theorem~4.4.1]:

\begin{longlist}[(I)]
\item[(I)] For $j=1,\ldots, J$, let $\psi_j = \vfi_{2j-1} + \ii\vfi
_{2j}$, where $\vfi_{1}, \ldots, \vfi_{2J} \in L^2( [0,1], \mathbb{R} )$, and $\b{Y}_t =
(Y_t(1), \ldots, Y_t(2J)) \in\mathbb{R}^{2J}$ be the vector time
series with coordinates $Y_t(l) = \sc{X_t, \vfi_l}.$
Then $\udf{\b{Y}}{\omega_{j,T}}^{(T)} \stackrel{d}{\rightarrow} \udf{\b
{Y}}{\omega_{j}}$, where $ \{ \udf{\b{Y}}{\omega_{j}}  \}$ are
independent mean zero complex Gaussian random vectors with covariance
matrix $\b{F}_{\omega_j},$ $(\b{F}_{\omega_j})_{sl}:= F_{\omega
_j}(s,l) = \sc{\specdO{\omega_j} \vfi_l, \vfi_s}.$
\end{longlist}

For the case $\mu\neq0$, we only need to consider $j=1,2$ since $\udf
{(X - \mu)}{\omega_{j,T}}^{(T)} = \udf{X}{\omega_{j,T}}^{(T)}$ for
$j=3, \ldots, J$. We need to show that
%
\begin{equation}
\label{eqthmudf-asymptotics-mu1} \udf{X} {\omega_{1}}^{(T)}
- \sqrt{\frac{T}{2\pi}} \mu= (2\pi T)^{-1/2} \sum
_{t=0}^{T-1} (X_{t} - \mu) \stackrel{d} {
\rightarrow} \udf{X} {0},
\end{equation}
and also that
%
\begin{equation}
\label{eqthmudf-asymptotics-mu2} \udf{X} {\omega_{2}}^{(T)}
= (2\pi T)^{-1/2} \sum_{t=0}^{T-1}
(-1)^{t} X_{t} \stackrel{d} {\rightarrow} \udf{X} {\pi}.
\end{equation}

The weak convergence in \eqref{eqthmudf-asymptotics-mu1} follows
immediately from the case $\mu= 0$. For~\eqref
{eqthmudf-asymptotics-mu2}, notice that
\[
\udf{X} {\omega_{2}}^{(T)} = (2\pi T)^{-1/2} \sum
_{t=0}^{T-1} (-1)^{t}
(X_{t} - \mu) + \mu (2\pi T)^{-1/2} \sum
_{t=0}^{T-1} (-1)^{t}.
\]
The first summand is the discrete Fourier transform of a zero mean
random process, and converges\vspace*{1pt} to $\udf{X}{\omega_{2}}$. The second
summand is deterministic and bounded by $\hnorm{\mu} (2\pi T)^{-1/2}$,
which tends to zero. Finally, the continuous mapping theorem for metric
spaces [\citet{Pollard1984}] yields \eqref{eqthmudf-asymptotics-mu2}.
\end{pf*}

The theorem has important consequences for the statistical analysis of
a functional time series. It essentially allows us to transform a
collection of weakly dependent functional data of an unknown
distribution,\vspace*{1pt} to a collection of approximately independent and Gaussian
functional data. In particular, let $\{\omega_{j,T}\}_{j=1}^{J}$ be\vspace*{-2pt} $J$
sequences (in $T$) of frequencies, such that $\omega_{j,T}\stackrel
{T\rightarrow\infty}{\longrightarrow}\omega\neq0$, for all $1\leq
j\leq J$. Then, provided $T$ is large enough, $\{\udf{X}{\omega
_{j,T}}^{(T)}\}_{j=1}^{J}$ is a collection of $J$\vspace*{1pt} approximately i.i.d.
mean zero complex Gaussian random functions with covariance operator
$\mathscr{F}_{\omega}$. The size $J$ of the sample is not allowed to
grow with~$T$, however. From a practical point of view, it can be
chosen to be large, provided that the $\omega_{j,T}$ are not too far
from $\omega$. We will make heavy use of this result in order to
construct consistent and asymptotically Gaussian estimators of the
spectral density operator by means of the \emph{periodogram kernel},
defined in the next section.

We also remark that the weak convergence relation in equation \eqref
{eqdftclt} can be re-expressed to trivially yield the corollary:
%
\begin{cor}[(Central limit theorem for {cumulant} mixing functional
series)] \label{corCLT}
Let $\{X_t\}_{t=0}^{T}$ be a strictly stationary sequence of random
elements of $ L^2( [0,1], \mathbb{R} )$ of length $T$ satisfying conditions \textup{(i)} and \textup{(ii)} of
Theorem~\ref{thmudf-asymptotics}. Then
\[
\sqrt{T} \Biggl(\frac{1}{T}\sum_{t=0}^{T-1}
X_{t}(\tau)-\mu(\tau) \Biggr)\stackrel{d} {\longrightarrow}\mathcal{N}
\biggl(0,\sum_{t\in\mathbb
{Z}}\mathscr{R}_t
\biggr).
\]
\end{cor}
This provides one of the first instances of central limit
theorems for functional series under no structural modelling
assumptions beyond weak dependence. To our knowledge, the only other
similar result is given in recent work by \citet{Horvath11}, who
obtain the same limit under different weak dependence conditions,
namely $L^p$-$m$-approximability.

The covariance operator of the limiting Gaussian measure is the
functional analogue of the long-run covariance matrix from multivariate
time series. We will revisit this operator in Section \ref
{secSpecdensityestimation}, where we will derive a related central
limit theorem.

\subsection{The periodogram kernel and its properties}\label{secperiodogram}

The covariance structure of the weak limit of the fDFT given in Theorem
\ref{thmudf-asymptotics} motivates the consideration of the empirical
covariance of the functional DFT as a basis for the estimation of the
spectral density operator. Thus, as with the multivariate case, we are
led to consider tensor products of the fDFT leading to the notion of a
\emph{periodogram kernel}. Define the \emph{periodogram kernel} as
\[
p_{\omega}^{(T)}(\tau, \sigma) = \bigl[\udf{X} {
\omega}^{(T)}(\tau) \bigr] \bigl[\udf{X} {\omega}^{(T)}(\sigma)
\bigr]^{\dag}=\udf{X} {\omega }^{(T)}(\tau) \udf{X} {-
\omega}^{(T)}(\sigma).
\]
If we slightly abuse notation and also write $\hnorm{\cdot}_2$ for the
norm in $L^{2}([0,1]^{2}, \mathbb C)$, we have
$\hnorm{p_{\omega}^{(T)}}_2 = \hnorm{\udf{X}{\omega}^{(T)}}^{2}_2,$
and hence
$\ee\hnorm{p_{\omega}^{(T)}}^{l}_2 < \infty$, if $\ee\hnorm
{X_t}_2^{2l} < \infty$.
The expectation of the periodogram kernel is thus well defined, and,
letting $a_{T} = \sum_{t=-T}^{T} e^{-\ii\omega t} \covK_{t}$,
Lemma~F.3 yields
$\ee p_{\omega}^{(T)} = T^{-1} (a_{0} + a_{1} + \cdots+ a_{T-1}).$
That is, the expectation of the periodogram kernel is a Ces\`aro-sum of
the partial sums of the series defining the spectral density kernel.
Therefore, in order to probe the properties of the periodogram kernel,
we can make use of the Fej\'er kernel
\[
F_{T}(\omega) = \frac{1}{T} \biggl( \frac{\sin(T\omega/ 2)}{\sin(\omega
/ 2)}
\biggr)^{2}=\frac{1}{ T} \bigl|\Delta^{(T)}(
\omega)\bigr|^{2}.
\]

It will thus be useful to recall some properties of $F_{T}$: $\int_{-\pi
}^{\pi} F_{T} = 2\pi, F_{T}(0) = T, F_{T}(\omega) \sim O(T)$
uniformly in $\omega$, and $F_{T}(2\pi s / T) = 0$ for $s$ an integer
with $s \not\equiv0 \mod T.$ This last property will be used often.
We will also be making use of the following cumulant mixing condition,
defined for fixed $l \geq0$ and $k =2,3,\ldots.$
\renewcommand{\thecondition}{$\mathrm{C}(l,k)$}
\begin{condition}
For each $j=1, \ldots, k-1,$
\[
\sum_{t_{1}, \ldots, t_{k-1} = -\infty}^{\infty} \bigl(1 +
|t_{j}|^{l} \bigr) \bigl\llVert \operatorname{cum}
(X_{t_{1}}, \ldots, X_{t_{k-1}}, X_{0} ) \bigr\rrVert
_{2} < \infty.
\]
\end{condition}
With this definition in place, we may determine the exact mean of the
periodogram kernel:
%
\begin{prop} \label{propperiodogram-exact-mean}
Assuming that $\mathrm{C}(0,2)$ holds true, we have, for each $\omega
\in\mathbb R,$
\[
\mathbb{E} \bigl[p_{\omega}^{(T)}(\tau, \sigma) \bigr] =
\frac{1}{2\pi} \int_{-\pi}^{\pi
} F_{T}(
\omega- \alpha)\specdK{\alpha}(\tau, \sigma) \,d\alpha+ \frac{1}{2\pi}\mu(\tau )\mu(\sigma)
F_{T}(\omega) \qquad\mbox{in } L^{2}.
\]
In particular, if $\omega= 2\pi s /T,$ with $s$ an integer such that
$s \not\equiv0 \mod T,$
\[
\mathbb{E} \bigl[ p_{\omega}^{(T)}(\tau, \sigma) \bigr] =
\frac{1}{2\pi} \int_{-\pi}^{\pi
} F_{T}(
\omega- \alpha)\specdK{\alpha}(\tau, \sigma) \,d\alpha\qquad \mbox{in }
L^{2}.
\]
\end{prop}

\begin{pf}
See the supplementary material [\citet{supp}], Proposition~C.1.\vadjust{\goodbreak}
\end{pf}

In particular, the periodogram kernel is asymptotically unbiased:

\begin{prop} \label{propapprox-expectation-periodogram}
Let $s$ be an integer with $s \not\equiv0 \mod T$. Then, we have
\[
\mathbb{E} \bigl[ p_{2\pi s /T}^{(T)}(\tau, \sigma) \bigr] =
\specdK{2\pi s/T}(\tau, \sigma) + \varepsilon_T  \qquad\mbox{in }
L^{2}.
\]
The error term $\varepsilon_T$ is $O(T^{-1})$ under $\mathrm{C}(0,2)$ and
$o(1)$ under $\mathrm{C}(1,2)$. In either case, the error term is uniform
in integers $s \not\equiv0 \mod T$.
\end{prop}

\begin{pf}
Since $s \not\equiv0 \mod T$,
\[
\mathbb{E} \bigl[p_{2\pi s/T}^{(T)}(\tau, \sigma) \bigr] =
\operatorname{cum} \bigl( \udf{X} {2\pi s/T}^{(T)}(\tau), \udf{X} {-2\pi
s/T}^{(T)}(\sigma) \bigr) = \specdK{2\pi s/T}(\tau, \sigma) +
\vep_t,
\]
and the result follows from Theorem~B.2 of
the supplementary material [\citet{supp}].
\end{pf}

Having established the mean structure of the periodogram, we turn to
the determination of its covariance structure.
%
\begin{thmm} \label{thmcov-raw-periodogram}
Assume $\omega_{1}$ and $\omega_{2}$ are of the form $2 \pi s(T) / T$,
where $s(T)$ is an integer, $s(T) \not\equiv0 \mod T$.
We have
\begin{eqnarray}
\cov \bigl( p_{\omega_{1}}^{(T)}(\tau_{1},
\sigma_{1}), p_{\omega
_{2}}^{(T)}(\tau_{2},
\sigma_{2}) \bigr) &=& \eta(\omega_{1} - \omega
_{2}) \specdK{\omega_{1}}(\tau_{1},
\tau_{2}) \specdK{- \omega _{1}}(\sigma_{1},
\sigma_{2})
\nonumber\\
&&{}+ \eta(\omega_{1} + \omega_{2}) \specdK{\omega
_{1}}(\tau_{1}, \sigma_{2}) \specdK{-
\omega_{1}}(\sigma_{1}, \tau_{2}) +
\vep_{T} \nonumber\\
\eqntext{\mbox{in } L^{2},}
\end{eqnarray}
where the function $\eta(x)$ equals one if $x \in2\pi\mathbb Z$, and
zero otherwise. The error term $\vep_{T}$ is $o(1)$ under $\mathrm{C}(0,2)$ and $\mathrm{C}(0,4)$; $\vep_{T} \sim O(T^{-1})$ under
$\mathrm{C}(1,2)$ and $\mathrm{C}(1,4)$. In each case, the error term is uniform
in $\omega_{1}, \omega_{2}$ [of the form $2\pi s(T)/T$ with $s(T) \not
\equiv0 \mod T$].

\end{thmm}

\begin{pf}
See the supplementary material [\citet{supp}], Theorem~C.2.
\end{pf}

\section{Estimation of the spectral density operator}\label
{secSpecdensityestimation}

The results in the previous section show that the asymptotic covariance
of the periodogram is not zero, and hence, as in the multivariate case,
the periodogram kernel itself is not a consistent estimator of the
spectral density. In this section, we define a consistent estimator,
obtained by convolving the periodogram kernel with an appropriate
weight function $W$. To this aim, let $W(x)$ be a real function defined
on $\mathbb R$ such that:
\begin{longlist}[(1)]
\item[(1)] $W$ is positive, even, and bounded in variation;
\item[(2)] $W(x) = 0$ if $|x| \geq1$;
\item[(3)] $\int_{-\infty}^{\infty} W(x)\ud x = 1;$
\item[(4)] $\int_{-\infty}^{\infty} W(x)^{2} \ud x < \infty.$\vadjust{\goodbreak}
\end{longlist}
The assumption of a compact support is not necessary, but will simplify
proofs. For a bandwidth $B_{T} > 0$, write
%
\begin{equation}
\label{eqkernel-of-order-T} W^{(T)}(x)= \sum
_{j \in\mathbb Z} \frac{1}{B_{T}}W \biggl( \frac{x + 2\pi j}{B_{T}} \biggr).
\end{equation}
Some properties of $W$ can be found in the supplementary material [\citet{supp}]. We define
the \emph{spectral density estimator} $\specdK{\omega}^{(T)}$ of
$\specdK{\omega}$ at frequency~$\omega$ as the weighted average of the
periodogram evaluated at frequencies of the form $\{2\pi s/T\}
_{s=1}^{T-1}$, with weight function $W^{(T)}$,
\[
\specdK{\omega}^{(T)}(\tau, \sigma) = \frac{2\pi}{T} \sum
_{s=1}^{T-1} W^{(T)} \biggl( \omega-
\frac{2\pi s}{T} \biggr) p_{{2\pi s}/{T}
}^{(T)}(\tau, \sigma).
\]
A consequence of the assumption of compact support worth mentioning is
that, in fact, at most $O(T B_T)$ summands of this expression are nonzero.
We will show in this section that, under appropriate conditions on the
asymptotic behavior of $B_T$, this estimator retains the property of
asymptotic unbiasedness that the periodogram enjoys. We will determine
the behaviour of its asymptotic covariance structure and establish
consistency in mean square (with respect to the Hilbert--Schmidt norm).
Finally, we will determine the asymptotic law of the estimator.

Concerning the mean of the spectral density estimator, we have:

\begin{prop} \label{propexpectation-specdk-estimator}
Under $\mathrm{C}(1,2)$,
if $B_{T}\rightarrow0$ and $B_{T}T \rightarrow\infty$ as $T
\rightarrow\infty$, then
\[
\ee\specdK{\omega}^{(T)}(\tau, \sigma) = \int_{\mathbb R}
W(x) \specdK {\omega- x B_{T}}(\tau, \sigma) \ud x + O
\bigl(B_{T}^{-1}T^{-1} \bigr),
\]
where the equality holds in $L^{2}$, and the error terms are uniform in
$ \omega.$
\end{prop}

\begin{pf}
See the supplementary material [\citet{supp}], Proposition~D.1.
\end{pf}

Concerning the covariance of the spectral density estimator,
we have:
%
\begin{thmm} \label{thmcov-specdk-estimator}
Under $\mathrm{C}(1,2)$ and $\mathrm{C}(1,4)$,
\begin{eqnarray*}
&& \cov \bigl( \specdK{\omega_{1}}^{(T)}(
\tau_{1}, \sigma_{1}), \specdK{\omega_{2}}^{(T)}(
\tau_{2}, \sigma_{2}) \bigr)
\\
&&\qquad= \frac{2\pi}{T} \int_{-\pi}^{\pi} \bigl\{
W^{(T)}(\omega_{1} - \alpha ) W^{(T)}(
\omega_{2} - \alpha) \specdK{\alpha}(\tau_{1},
\tau_{2}) \specdK{-\alpha}(\sigma_{1}, \sigma_{2})
\\
&&\hspace*{40pt}\qquad\quad{}+ W^{(T)}(\omega_{1} - \alpha) W^{(T)}(
\omega_{2} + \alpha) \specdK{\alpha}(\tau_{1},
\sigma_{2}) \specdK{-\alpha}(\sigma _{1},
\tau_{2}) \bigr\}\ud\alpha
\\
&&\qquad\quad{} + O \bigl(B_{T}^{-2} T^{-2} \bigr) + O
\bigl(T^{-1} \bigr),
\end{eqnarray*}
where the equality holds in $L^{2}$, and the error terms are uniform in
$\omega.$

\end{thmm}

\begin{pf}
See the supplementary material [\citet{supp}], Theorem~D.2.
\end{pf}

Noting that $\| W^{(T)} \|_{\infty} = O(B_{T}^{-1})$ and $\| \specdK
{\cdot} \|_{\infty} = O(1)$, a direct consequence of the last result is
the following approximation of the asymptotic covariance of the
spectral density estimator:
%
\begin{cor} \label{corcov-unif-bound}
Under $\mathrm{C}(1,2)$ and $\mathrm{C}(1,4)$,
\[
\cov \bigl( \specdK{\omega_{1}}^{(T)}(\tau_{1},
\sigma_{1}), \specdK {\omega_{2}}^{(T)}(
\tau_{2}, \sigma_{2}) \bigr) = O \bigl( B_{T}^{-2}
T^{-1} \bigr),
\]
where the equality holds in $L^{2}$, uniformly in the $\omega$'s.
\end{cor}

This bound is not sharp. A better bound is given in the next
statement, which, however, is not uniform in $\omega.$
%
\begin{prop} \label{propcov-sharp-bound}

Assume conditions $\mathrm{C}(1,2)$, $\mathrm{C}(1,4)$, and that $B_{T}
\rightarrow0$ as $T \rightarrow\infty$ with $B_{T} T \rightarrow
\infty.$
Then
\begin{eqnarray*}
&&\lim_{T \rightarrow\infty} B_{T} T \cov \bigl( \specdK {
\omega_{1}}^{(T)}(\tau_{1}, \sigma_{1}),
\specdK{\omega_{2}}^{(T)}(\tau _{2},
\sigma_{2}) \bigr) \\
&&\qquad= 2 \pi\int_{\mathbb R} W(
\alpha)^{2} \ud \alpha\bigl\{ \eta(\omega_{1} - \omega_{2}) \specdK{
\omega _{1}}(\tau_{1}, \tau_{2}) \specdK{-
\omega_{1}}(\sigma_{1}, \sigma_{2})\\
&&\hspace*{110pt}{} + \eta(
\omega_{1} + \omega_{2}) \specdK{\omega_{1}}(
\tau_{1}, \sigma _{2}) \specdK{-\omega_{1}}(
\sigma_{1}, \tau_{2}) \bigr\}.
\end{eqnarray*}
The function $\eta(x)$ equals one if $x \in2\pi\mathbb Z$, and zero
otherwise. The convergence is in $L^{2}$ for any fixed $\omega_{1},
\omega_{2}.$ If $\omega_{1}, \omega_{2}$ depend on $T$, then the
convergence is in $L^{2}$, provided $(\omega_{1} \pm\omega_{2})$ are
at a distance of at least $2 B_{T}$ from any multiples of $2\pi$, if
not exactly a multiple of $2\pi$.
\end{prop}

\begin{pf}
Let $d(x,y)$ denote the distance in $\mathbb R / 2\pi\mathbb Z$. We
shall abuse notation and let $x,y$ stand for equivalence classes of
real numbers,
and also omit the $(\tau, \sigma)$'s, for the sake of clarity.
Theorem~\ref{thmcov-specdk-estimator} yields
%
\begin{eqnarray}
&&B_{T} T \cov \bigl( \specdK{\omega_{1}}^{(T)},
\specdK {\omega_{2}}^{(T)} \bigr)
\nonumber
\\
&&\qquad= 2\pi B_{T} \int_{-\pi}^{\pi}
W^{(T)}(\omega_{1} - \omega_{2} - \alpha)
W^{(T)}(\alpha) \specdK{\omega_{2} + \alpha} \specdK{-(\omega
_{2} + \alpha)} \,d\alpha\label{eqcorcov-sharp-bound-eq1}
\\
&&\qquad\quad{} + 2\pi B_{T} \int_{-\pi}^{\pi}
W^{(T)}(\omega_{1} + \omega_{2} - \alpha)
W^{(T)}( \alpha) \specdK{-(\omega_{2} - \alpha)} \specdK{
\omega_{2} -\alpha} \,d\alpha \label{eqcorcov-sharp-bound-eq2}
\\
&&\qquad\quad{} + O \bigl(B_{T}^{-1} T^{-1} \bigr) +
O(B_{T}).
\nonumber
\end{eqnarray}
We have employed a change of variables, the fact that $W^{(T)}$ is
even, and the fact that both $W^{(T)}$ and $\specdK{\cdot}$ are $2\pi
$-periodic. The error terms tend to zero as $B_{T} \rightarrow0$,
$TB_{T} \rightarrow\infty.$\vadjust{\goodbreak}

First we show that \eqref{eqcorcov-sharp-bound-eq1} tends to
%
\begin{equation}
\label{eqcorcov-sharp-bound-eq3} \eta(\omega_{1} -
\omega_{2}) \specdK{\omega_{1}}(\tau_{1},
\tau_{2}) \specdK{-\omega_{1}}(\sigma_{1},
\sigma_{2}) 2 \pi\int_{\mathbb R} W(\alpha)^{2}
\ud\alpha,
\end{equation}
in $L^{2}$, uniformly in all $\omega_{1} = \omega_{1,T}, \omega_{2} =
\omega_{2,T}$ such that $\omega_{1,T} \equiv\omega_{2,T}$ or $d(\omega
_{1,T} - \omega_{2,T}, 0) \geq2 B_{T}$ for large $T$.
If $d(\omega_{1} - \omega_{2}, 0) \geq2 B_{T}$, \eqref
{eqcorcov-sharp-bound-eq1} is exactly equal to zero.
If $\omega_{1} \equiv\omega_{2},$ we claim that \eqref
{eqcorcov-sharp-bound-eq1} tends to
%
\begin{equation}
\label{eqcorcov-sharp-bound-eq4} \specdK{\omega}(\tau_{1},
\tau_{2}) \specdK{-\omega}(\sigma_{1}, \sigma
_{2}) 2\pi\int_{\mathbb R} W(\alpha)^{2} \,d
\alpha.
\end{equation}
Notice that in this case, \eqref{eqcorcov-sharp-bound-eq1} can be
written as
$\int_{-\pi}^{\pi} K_{T}(\alpha) \specdK{\omega+ \alpha} \specdK
{-(\omega+ \alpha)} \,d\alpha\times \{ \int_{\mathbb R} W(\alpha
)^{2} \,d\alpha \},$
where $K_{T}(\alpha) = \frac{2\pi}{B_{T}}  [ W(\alpha/ B_{T})
]^{2}  \{ \int_{\mathbb R} W(\alpha)^{2} \,d\alpha \}^{-1}$
is an approximate identity on $[-\pi, \pi]$; see \citet{Edwards1967}, Section
3.2.
{Since the spectral density kernel is uniformly continuous with
respect to $\hnorm{\cdot}_{2}$ (see Proposition~\ref{lmaprop-specdk})}
Lemma~F.15 implies that \eqref
{eqcorcov-sharp-bound-eq1} tends indeed to \eqref
{eqcorcov-sharp-bound-eq4} uniformly in $\omega$ with respect to
$\hnorm{\cdot}_{2}$.
Hence \eqref{eqcorcov-sharp-bound-eq1} tends to \eqref
{eqcorcov-sharp-bound-eq3} in $\hnorm{\cdot}_{2}$, uniformly in
$\omega$'s satisfying
\[
\omega_{1,T} \equiv\omega_{2,T} \quad\mbox{or}\quad d(
\omega_{1,T} - \omega _{2,T}, 0) \geq2 B_{T}\qquad
\mbox{for large } T.
\]

Similarly, we may show that \eqref
{eqcorcov-sharp-bound-eq2} tends to
$\eta(\omega_{1} + \omega_{2}) \specdK{\omega_{1}}(\tau_{1}, \sigma
_{2}) \specdK{-\omega_{1}}(\sigma_{1}, \tau_{2}) \times 2 \pi\int_{\mathbb R}
W(\alpha)^{2} \ud\alpha,$
uniformly in $\omega$'s if $\omega_{1,T} \equiv-\omega_{2,T}$ or
$d(\omega_{1,T} + \omega_{2,T}) \geq2 B_{T}$ for large $T$.
Piecing these results together, we obtain the desired convergence,
provided for each $T$ large enough, either
$\omega_{1,T} - \omega_{2,T} \equiv0,$
$\omega_{1,T} + \omega_{2,T} \equiv0,$
or
\[
d(\omega_{1,T} - \omega_{2,T}, 0 ) \geq2 B_{T}
\quad\mbox{and}\quad d(\omega_{1,T} + \omega_{2,T}, 0) \geq2
B_{T}.
\]
\upqed\end{pf}

\begin{rmk}\label{rmksup-norm-cumulant-conditions}
In practice, functional data are assumed to be smooth in addition to
square-integrable. In such cases, one may hope to obtain stronger
results, for example with respect to uniform rather than $L^2$ norms.
Indeed, if the conditions {$\mathrm{C}(l,k)$} are replaced by the stronger conditions
\renewcommand{\thecondition}{$\mathrm{C}'(l,k)$}
\begin{condition}For each $j = 1,\ldots, k - 1$ 
\[
\sum_{t_{1}, \ldots, t_{k-1} \in\mathbb Z}
\bigl(1 + |t_{j}|^{l} \bigr) \bigl\Vert \operatorname{cum} (
X_{t_{1}}, \ldots, X_{t_{k-1}}, X_{0} ) \bigr\Vert_{\infty}
< \infty,
\]
\end{condition}
then the results of Propositions~\ref{propperiodogram-exact-mean}, \ref
{propapprox-expectation-periodogram}, Theorem~\ref
{thmcov-raw-periodogram}, Proposition~\ref
{propexpectation-specdk-estimator}, Theorem~\ref
{thmcov-specdk-estimator}, Corollary~\ref{corcov-unif-bound},
Proposition~\ref{propcov-sharp-bound},
and Lemma~B.1, Theorem~B.2
in the supplementary material [\citet{supp}] would hold in the supremum norm with respect to
$\tau, \sigma$.
\end{rmk}

Combining the results on the asymptotic bias and variance of the
spectral density operator, we may now derive the consistency in
integrated mean square of the induced estimator for the spectral
density operator.\vadjust{\goodbreak} Recall that $\specdO{\omega}$ is the integral
operator with kernel $\specdK{\omega},$ and, similarly let $\specdO
{\omega}^{(T)}$ be the operator with kernel $\specdK{\omega}^{(T)}$. We have:

\begin{thmm} \label{thmimse-convergence}
Provided assumptions $\mathrm{C}(1,2)$ and $\mathrm{C}(1,4)$ hold, $B_{T}
\rightarrow0$, $B_{T} T \rightarrow\infty,$
the spectral density operator estimator $\specdO{\omega}^{(T)}$ is
consistent in integrated mean square, that is,
\[
\operatorname{IMSE} \bigl(\specdO{}^{(T)} \bigr) = \int
_{-\pi}^{\pi} \ee\bigl \vert\!\bigl\vert\!\bigl\vert \specdO{
\omega}^{(T)} - \specdO{\omega} \bigr\vert\!\bigr\vert\!\bigr\vert^{2}_{2}
\ud\omega \rightarrow0, \qquad T \rightarrow\infty,
\]
where $\snorm{\cdot}_{2}$ is the Hilbert--Schmidt norm (the Schatten 2-norm).\vspace*{1pt}
More precisely,
$\operatorname{IMSE}(\specdO{}^{(T)}) = O(B_{T}^{2}) +
O(B_{T}^{-1}T^{-1})  \mbox{ as } T \rightarrow\infty.$
We also have pointwise mean square convergence for a fixed $\omega$:
\[
\ee \bigl\vert\!\bigl\vert\!\bigl\vert \specdO{\omega}^{(T)} - \specdO{\omega} \bigr\vert\!\bigr\vert\!\bigr\vert^{2}_{2} = \cases{ O \bigl(B_{T}^{2}
\bigr) + O \bigl(B_{T}^{-1}T^{-1} \bigr), & \quad$
\mbox{if } 0 < |\omega| < \pi,$ \vspace*{2pt}
\cr
O \bigl(B_{T}^{2}
\bigr) + O \bigl(B_{T}^{-2}T^{-1} \bigr), & \quad$
\mbox{if } \omega=0, \pm\pi$}
\]
as $T \rightarrow\infty.$

\end{thmm}

\begin{pf}
For an integral operator $K$ with a complex-valued\vspace*{1pt} kernel $k(\tau, \sigma
)$, we will denote by $\overline{K}$ the operator with kernel $\overline
{k(\tau, \sigma)}.$ Let $\snorm{\cdot}_{2}$ be the Hilbert--Schmidt norm.
Proposition F.21 yields
$\snorm{K}_2 = \snorm{ \overline{K}}_2.$ Further, notice that $\specdK
{-\omega}(\tau, \sigma) = \overline{\specdK{\omega} (\tau, \sigma)}$,
hence $\specdO{-\omega} = \overline{\specdO{\omega}}.$ Similarly,
$\specdO{-\omega}^{(T)} = \overline{\specdO{\omega}^{(T)}}.$ Thus, via
a change of variables, the IMSE of the spectral density estimator can
be written as
\begin{eqnarray*}
&& \int_{-\pi}^{\pi} \ee \bigl\vert\!\bigl\vert\!\bigl\vert
\specdO{\omega}^{(T)}- \specdO {\omega} \bigr\vert\!\bigr\vert\!\bigr\vert_{2}^{2}
\ud\omega= 2 \int_{0}^{\pi} \ee \bigl\vert\!\bigl\vert\!\bigl\vert
\specdO {\omega}^{(T)} - \specdO{\omega} \bigr\vert\!\bigr\vert\!\bigr\vert
_{2}^{2} \ud\omega
\\
&&\qquad= 2 \int_{0}^{\pi} \ee \bigl\vert\!\bigl\vert\!\bigl\vert \specdO{
\omega}^{(T)} - \ee\specdO {\omega}^{(T)} \bigr\vert\!\bigr\vert\!\bigr\vert
_{2}^{2} \ud\omega + 2\int_{0}^{\pi}
\bigl\vert\!\bigl\vert\!\bigl\vert \specdO{\omega} - \ee\specdO{\omega}^{(T)} \bigr\vert\!\bigr\vert\!\bigr\vert _{2}^{2} \ud\omega,
\end{eqnarray*}
which is essentially the usual bias/variance decomposition of the mean
square error. Initially, we focus on the variance term. Lemma~F.22 yields
\[
\int_{0}^{\pi} \ee \bigl\vert\!\bigl\vert\!\bigl\vert \specdO{
\omega}^{(T)} - \ee\specdO{\omega }^{(T)} \bigr\vert\!\bigr\vert\!\bigr\vert_{2}^{2} \ud\omega= \int_0^{\pi}
\iint_{[0,1]^{2}} \var \bigl( \specdK{\omega}^{(T)}(\tau, \sigma)
\bigr) \,d\tau \,d\sigma\ud\omega.
\]
Decomposing the outer integral into three terms, $\int_{0}^{\pi} = \int_{0}^{\pi B_{T}} + \int_{\pi B_{T}}^{\pi-B_{T}} + \int_{\pi-B_{T}}^{\pi
},$ we can use Corollary~\ref{corcov-unif-bound} for the first and
last summands, and Proposition~\ref{propcov-sharp-bound} for the
second summand to obtain
$ \int_{0}^{\pi} \ee\snorm{ \specdO{\omega}^{(T)} - \ee\specdO{\omega
}^{(T)} }_{2}^{2} \ud\omega= O(B_{T}^{-1} T^{-1}).$ Turning to the
squared bias, Proposition~\ref{propexpectation-specdk-estimator} yields
\begin{eqnarray*}
&& \int_{0}^{\pi} \bigl\vert\!\bigl\vert\!\bigl\vert
\specdO{\omega} - \ee\specdO{\omega }^{(T)} \bigr\vert\!\bigl\vert\!\bigr\vert
_{2}^{2} \ud\omega
\\
&&\qquad\leq3\int_{0}^{\pi} \biggl\vert\!\biggl\vert\!\biggl\vert \biggl\{
\int_{\mathbb R} W(x) \specdK {\omega- x B_{T}} \ud x -
\specdK{\omega} \biggr\} \biggr\vert\!\biggr\vert\!\biggr\vert _{2}^{2} \ud
\omega + O \bigl(T^{-2} \bigr) + O \bigl(B_{T}^{-2}T^{-2}
\bigr),
\end{eqnarray*}
where we have used Jensen's inequality and where
$ \{ \int_{\mathbb R} W(x) \specdK{\omega- x B_{T}} \ud x -
\specdK{\omega}  \}$
denotes the operator with kernel
$ \int_{\mathbb R} W(x) \specdK{\omega- x B_{T}}(\tau, \sigma) \ud x
- \specdK{\omega}(\tau, \sigma).$
Lem\-ma~F.4 implies that this
difference is of order $O(B_{T})$, uniformly in $\omega.$ Hence,
\[
3 \int_{0}^{\pi} \biggl\vert\!\biggl\vert\!\biggl\vert \biggl\{ \int
_{\mathbb R} W(x) \specdK{\omega- x B_{T}} \ud x -
\specdK{\omega} \biggr\} \biggr\vert\!\biggr\vert\!\biggr\vert _{2}^{2} \ud
\omega\leq O \bigl(B_{T}^{2} \bigr).
\]
In summary, we have
\[
\int_{-\pi}^{\pi} \ee\bigl\vert\!\bigl\vert\!\bigl\vert\specdO{
\omega}^{(T)} - \specdO{\omega }\bigr\vert\!\bigr\vert\!\bigr\vert^{2} \ud\omega\leq O
\bigl(B_{T}^{2} \bigr) + O \bigl(B_{T}^{-1}T^{-1}
\bigr).
\]
The spectral density estimator $\specdO{\cdot}^{(T)}$ is therefore
consistent in integrated mean square if $B_{T} \rightarrow0$ and
$B_{T}T \rightarrow\infty$ as $T\rightarrow\infty$.

A careful examination of the proof reveals that the
pointwise statement of the theorem follows by a directly analogous argument.
\end{pf}

Finally, if we include some higher-order cumulant mixing
conditions, we may obtain the asymptotic distribution of our estimator
as being Gaussian.

\begin{thmm}\label{thmCLTspecdensity}
Assume that {$\ee\hnorm{X_0}^k < \infty$ for all $k \geq2$ and}:
\begin{longlist}[(iii)]
\item[(i)]
$\sum_{t_{1}, \ldots, t_{k-1} = -\infty}^{\infty} \llVert  \cum
{X_{t_{1}}, \ldots, X_{t_{k-1}}, X_{0}} \rrVert _{2} < \infty, $
for all $k\geq2; $
\item[(i$'$)]
$\sum_{t_{1}, \ldots, t_{k-1} = -\infty}^{\infty} (1+|t_j|)\llVert  \cum
{X_{t_{1}}, \ldots, X_{t_{k-1}}, X_{0}} \rrVert _{2} < \infty,$ for
$k\in\{2,4\}$ and $j<k$;
\item[(ii)] $\sum_{t \in\mathbb Z}(1 + |t|) \snorm{\covO_{t}}_{1} <
\infty;$
\item[(iii)] $ \sum_{t_{1}, t_{2}, t_{3} \in\mathbb Z} \snorm{ \covO
_{t_{1}, t_{2}, t_{3}} }_{1} < \infty.$
\end{longlist}

Then, for any frequencies $\omega_{1}, \ldots, \omega_{J}$,
with $J < \infty$,
\[
\sqrt{B_{T}T} \bigl( \specdK{\omega_{j}}^{(T)} -
\ee \specdK{\omega _{j}}^{(T)} \bigr) \stackrel{d} {
\longrightarrow} \specdKE{\omega_{j}},\qquad j=1,\ldots, J,
\]
where $\specdKE{\omega_{j}}, j=1, \ldots, J,$ are jointly mean zero
complex Gaussian elements in $L^{2} ( [0,1]^{2}, \mathbb C
)$, with covariance kernel
\begin{eqnarray*}
&&\cov \bigl( \specdKE{\omega_{i}}(\tau_{1},
\sigma_{1}), \specdKE{\omega_{j}}(\tau_{2},
\sigma_{2}) \bigr)\\
&&\qquad= 2\pi\int_{\mathbb
R} W(
\alpha)^{2} \,d\alpha \bigl\{ \eta(\omega_{i} - \omega_{j}) \specdK{
\omega_{i}}(\tau_{1}, \tau_{2}) \specdK{-
\omega_{i}}(\sigma _{1}, \sigma_{2})\\
&&\hspace*{88pt}\qquad{} + \eta(
\omega_{i} + \omega_{j}) \specdK{\omega_{i}}(
\tau_{1}, \sigma _{2}) \specdK{-\omega_{i}}(
\sigma_{1}, \tau_{2}) \bigr\}.
\end{eqnarray*}
In particular, we see that $\specdKE{\omega_{i}}$ and $\specdKE{\omega
_{j}}$ are independent if $\omega_{i} \pm\omega_{j} \not\equiv0 \mod
2\pi,$ and $\specdKE{\omega}$ is real Gaussian if $\omega\equiv0 \mod
\pi.$
\end{thmm}

\begin{pf}
Let $(\vfi_{n})$ be a basis of $ L^2( [0,1], \mathbb{R} )$. Then $\{\vfi_{m_{1}} \otimes
\cdots\otimes\break\vfi_{m_{k}}\}_{m_{1}, \ldots, m_{k} \geq1}$ is a basis
of the complex Hilbert space $L^{2}([0,1]^{k}, \mathbb C)$ [e.g., \citet
{Kadison1997}], where $\vfi_{m_{1}} \otimes\cdots\otimes\vfi
_{m_{k}}(\tau_1, \ldots, \tau_k) = \vfi_{m_{1}}(\tau_1) \times\cdots
\times \vfi_{m_{k}}(\tau_k).$\vadjust{\goodbreak}
We denote by $\specdM{\omega}(m,n)$ the $(m,n)$th coordinate of the
spectral density matrix and, more generally, we define the $k$th-order
cumulant spectra array $\Phi_{\omega_{1}, \ldots, \omega_{k-1}}$ by
\[
\Phi_{\omega_{1}, \ldots, \omega_{k-1}}(m_{1}, \ldots, m_{k}) = \int
_{[0,1]^{k}} \specdK{\omega_{1}, \ldots,
\omega_{k-1}}(\tau_{1}, \ldots, \tau_{k})
\vfi_{m_{1}}(\tau_{1}) \cdots\vfi_{m_{k}}(
\tau_{k}) \,d\boldsymbol{\tau}.
\]
In other words, the $k$th-order cumulant spectra array is the scalar
product in $L^{2}([0,1]^{k}, \mathbb C)$ between $\specdK{\omega_{1},
\ldots, \omega_{k-1}}$ and $\vfi_{m_{1}} \otimes\cdots\otimes\vfi_{m_{k}}$.

We also define the periodogram matrix $P_{\omega}^{(T)}(m,n)$ and the
estimator of the spectral density matrix $\specdM{\omega}^{(T)}(m,n)$
as the scalar product of $p_{\omega}^{(T)}$, respectively, $\specdK
{\omega}^{(T)}$, with $\vfi_{m}\otimes\vfi_{n}.$ Notice that
\[
\specdM{\omega}^{(T)}(m,n) = \frac{2\pi}{T} \sum
_{s=1}^{T-1} W^{(T)} (\omega- 2\pi s/T)
P^{(T)}_{2\pi s/T}(m,n),
\]
where $W^{(T)}$ is defined in \eqref{eqkernel-of-order-T}. The major
steps in the proof are the following two inequalities, along with Lemma
\ref{lmatech-criterion-unif-tightness}:
\begin{longlist}[(II)]
\item[(I)] For $T$ large enough,
%
\begin{eqnarray}\label{eqspecdM-explicit-var}
&&T B_{T} \var \bigl( \specdM{\omega}^{(T)}(m,n) \bigr)
\nonumber
\\
&&\qquad  \leq
K \bigl[ \operatorname{sc}_{0}(m,n,m,n) + \operatorname{sc}_{1}(m,m)\operatorname{sc}_{1}(n,n)
\\
&&\hspace*{20pt}\qquad\quad{}+ \operatorname{sc}_{0}(m,m)\operatorname{sc}_{0}(n,n) + \operatorname{sc}_{0}(m,n)^{2}
\bigr],\nonumber
\end{eqnarray}
where ``sc'' stands for ``summed cumulant,'' in particular
\begin{eqnarray*}
\operatorname{sc}_{0}(m_{1}, \ldots, m_{k}) &=& \sum
_{t_{1}, \ldots, t_{k-1} \in\mathbb
Z} \bigl\llvert \operatorname{cum} \bigl(
\xi_{t_{1}}(m_{1}), \ldots, \xi_{t_{k-1}}(m_{k-1}),
\xi _{0}(m_{k}) \bigr) \bigr\rrvert,
\\
\operatorname{sc}_{1}(m_{1}, \ldots, m_{k}) &=& \sum
_{j=1}^{k-1}\sum_{t_{1}, \ldots,
t_{k-1} \in\mathbb Z}
|t_{j}| \bigl\llvert \operatorname{cum} \bigl( \xi_{t_{1}}(m_{1}),
\ldots, \xi_{t_{k-1}}(m_{k-1}), \xi_{0}(m_{k})
\bigr) \bigr\rrvert,
\end{eqnarray*}
and $K = 8 \hnorm{W}_{\infty}^{2}$ is a constant.

\item[(II)] We have the following bound:
%
\begin{eqnarray}\label{eqbound-sum-var-specdKE}
&&\sum_{m,n \geq1}T B_{T} \var \bigl( \specdM{
\omega}^{(T)}(m,n) \bigr)
\nonumber
\\[-8pt]
\\[-8pt]
\nonumber
&&\qquad\leq C \biggl[ \biggl( \sum
_{t\in\mathbb Z}\bigl(1 + |t|\bigr) \vert\!\vert\!\vert \covO _{t} \vert\!\vert\!\vert _{1} \biggr)^{2} + \sum_{t_{1}, t_{2}, t_{3} \in\mathbb
Z} \vert\!\vert\!\vert
\covO_{t_{1}, t_{2}, t_{3}} \vert\!\vert\!\vert _{1} \biggr]
\end{eqnarray}
for some fixed constant $C > 0.$ Here, $\covO_{t_{1}, t_{2}, t_{3}}$ is
the operator on $L^{2} ([0,1]^{2}, \mathbb R  )$ with kernel
$\covK_{t_{1}, t_{2}, t_{3}}  ( (\tau_{1}, \tau_{2}), (\tau_{3},
\tau_{4})  ) = \cum{X_{t_{1}}, X_{t_{2}}, X_{t_{3}}, X_{0} }(\tau
_{1}, \tau_{2}, \tau_{3}, \tau_{4}).$\break
That is,
$\covO_{t_{1}, t_{2}, t_{3}} f (\tau_{1}, \tau_{2}) = \iint
_{[0,1]^{2}} \covK_{t_{1}, t_{2}, t_{3}}  ( (\tau_{1}, \tau_{2}),
(\tau_{3}, \tau_{4})  ) f(\tau_{3}, \tau_{4}) \,d\tau_{3} \,d\tau_{4}$
for $f \in L^{2}([0,1]^{2}, \mathbb R).$

\end{longlist}

First we concentrate on establishing \eqref{eqspecdM-explicit-var}.
Recall that
\begin{eqnarray*}
\var \bigl( \specdM{\omega}^{(T)}(m,n) \bigr) &=& (2\pi/
T)^{2} \sum_{s,l = 1}^{T-1}
W^{(T)}(\omega- 2\pi s/T) W^{(T)}(\omega- 2\pi l/T)
\\
&&\hspace*{60pt}{}\times \cov \bigl( P_{{2\pi s}/{T}}^{(T)}(m,n), P_{{2\pi l}/{T}}^{(T)}(m,n)
\bigr).
\nonumber
\end{eqnarray*}
We need to find an explicit bound on the error terms of Lemma~B.1, Theorem~B.2 in the supplementary material [\citet{supp}], and Theorem~\ref{thmcov-raw-periodogram}. An examination
of the proof of Lemma~B.1 in the supplementary material [\citet{supp}] yields
\begin{eqnarray*}
&&\specdM{\omega_{1}, \ldots, \omega_{k-1}}(m_{1},
\ldots, m_{k})\\
&&\qquad= (2\pi)^{-(k-1)} \sum
_{t_{1}, \ldots, t_{k-1} =
-(T-1)}^{T-1} \exp \Biggl(- \ii\sum
_{j=1}^{k-1} \omega _{j}t_{j}
\Biggr)\\
 &&\hspace*{118pt}\qquad\quad{}\times\operatorname{cum} \bigl( \xi_{t_{1}}(m_{1}),
\ldots, \xi _{t_{k-1}}(m_{k-1}), \xi_{0}(m_{k})
\bigr) \\
&&\hspace*{118pt}\qquad\quad{}+ \vep^{\mathrm{(B.1)}}_{T}(m_{1}, \ldots,
m_{k}),
\end{eqnarray*}
and
$ \llvert  \vep^{\mathrm{(B.1)}}_{T}(m_{1}, \ldots, m_{k})
\rrvert  \leq(2\pi)^{-(k-1)} (k-1)   \operatorname{sc}_{0}(m_{1}, \ldots, m_{k}).$
We have used the notation $\vep_{T}^{\mathrm{(B.1)}}(m_{1}, \ldots, m_{k})$ to denote the error term
of Lemma~B.1, and we shall do likewise for the
error term in Theorem~B.2 in the supplementary material [\citet{supp}],
\begin{eqnarray*}
&&T^{k/2} \operatorname{cum} \bigl( \udf{\xi} {
\omega_{1}}^{(T)}(m_{1}), \ldots, \udf {\xi} {
\omega_{k}}^{(T)}(m_{k}) \bigr)
\\
&&\qquad= (2\pi)^{k/2 -1} \Delta^{(T)} \Biggl( \sum
_{j=1}^{k} \omega_{j} \Biggr) \specdM{
\omega_{1}, \ldots, \omega_{k-1}}(m_{1}, \ldots,
m_{k})
\\
&&\qquad\quad{} + \vep^{\mathrm{(B.2)}}_{T} \Biggl( \sum
_{j=1}^{k} \omega_{j}; m_{1},
\ldots, m_{k} \Biggr),
\end{eqnarray*}
where
\begin{eqnarray*}
&& \bigl\llvert \vep^{\mathrm{(B.2)}}_{T} (\omega;
m_{1}, \ldots, m_{k} ) \bigr\rrvert
\\[-1pt]
&&\qquad\leq2 (2\pi)^{-k/2} \sum_{t_{1}, \ldots,
t_{k-1}=-(T-1)}^{T-1}
\bigl( |t_{1}| + \cdots+ |t_{k-1}| \bigr) \\[-1pt]
&&\hspace*{128pt}\qquad{}\times\bigl\llvert\operatorname{cum} \bigl( \xi_{t_{1}}(m_{1}), \ldots, \xi
_{t_{k-1}}(m_{k-1}), \xi_{0}(m_{k}) \bigr)
\bigr\rrvert
\\[-1pt]
&&\hspace*{120pt}\qquad\quad{} + (2\pi)^{k/2 -1} \Delta^{(T)} ( \omega ) \bigl\llvert
\vep_{T}^{\mathrm{(B.1)}}(m_{1}, \ldots, m_{k})
\bigr\rrvert
\\[-1pt]
&&\qquad\leq2 (2\pi)^{-k/2} \operatorname{sc}_{1}(m_{1}, \ldots,
m_{k}) + (2\pi )^{-k/2} (k-1) \Delta^{(T)} ( \omega )
\operatorname{sc}_{0}(m_{1}, \ldots, m_{k}).
\end{eqnarray*}
A less sharp bound (but independent of the frequency) will also be useful,
\[
\bigl\llvert \vep^{\mathrm{(B.2)}}_{T} ( \cdot; m_{1},
\ldots, m_{k} ) \bigr\rrvert \leq3 (2\pi)^{-k/2} (k-1) T
\operatorname{sc}_{0}(m_{1}, \ldots, m_{k}).\vadjust{\goodbreak}
\]
We will also need a bound on the spectral density matrix,
$\llvert  \specdM{\omega_{1}, \ldots, \omega_{k-1}}(m_{1},
\ldots,\break
m_{k}) \rrvert  \leq(2\pi)^{-(k-1)}   \operatorname{sc}_{0}(m_{1}, \ldots, m_{k}).$

We now turn to Theorem~\ref{thmcov-raw-periodogram}: for $s,l=1,
\ldots, T-1,$
\begin{eqnarray*}
&&\cov \bigl( P_{{2\pi s}/{T}}^{(T)}(m,n), P_{{2\pi
l}/{T}}^{(T)}(m,n)
\bigr)
\\[-2pt]
&&\qquad= (2\pi/T) \specdM{{2\pi s}/{T}, -{2\pi s}/{T},
{2\pi
l}/{T}}(m,n,m,n) + T^{-2} \vep_{T}^{\mathrm{(B.2)}}(
\cdot;m,n,m,n)
\\[-2pt]
&&\qquad\quad{} + \delta_{s,l} \bigl[ \specdM{{2\pi
s}/{T}}(m,m)\specdM{-
{2\pi s}/{T}}(n,n) + \specdM{{2\pi
s}/{T}}(m,m) T^{-1}
\vep_{T}^{\mathrm{(B.2)}}(\cdot;n,n)
\\[-2pt]
&&\hspace*{156pt}\qquad\quad{} + \specdM{-{2\pi s}/{T}}(n,n) T^{-1} \vep_{T}^{\mathrm{(B.2)}}(
\cdot;m,m) \bigr]
\\[-2pt]
&&\qquad\quad{} + \delta_{s+l,T} \bigl[ \specdM{{2\pi
s}/{T}}(m,n)\specdM{-
{2\pi s}/{T}}(n,m) \\[-2pt]
&&\hspace*{55pt}\qquad{}+ \specdM{{2\pi
s}/{T}}(m,n) T^{-1}
\vep_{T}^{\mathrm{(B.2)}}(\cdot;n,m)
\\[-2pt]
&&\hspace*{45pt}\qquad\quad{} + \specdM{-{2\pi s}/{T}}(n,m) T^{-1} \vep_{T}^{\mathrm{(B.2)}}(
\cdot;m,n) \bigr]
\\[-2pt]
&&\qquad\quad{} + T^{-2} \biggl[ \vep_{T}^{\mathrm{(B.2)}} \biggl(
\frac{2\pi(s-l)}{T};m,m \biggr) \vep _{T}^{\mathrm{(B.2)}} \biggl(-
\frac{2\pi(s-l)}{T};n,n \biggr)
\\[-2pt]
&&\hspace*{35pt}\qquad\quad{} + \vep_{T}^{\mathrm{(B.2)}} \biggl( \frac{2\pi(s+l)}{T}; m,n \biggr)
\vep_{T}^{\mathrm{(B.2)}} \biggl( - \frac{2\pi(s+l)}{T}; n,m \biggr)
\biggr],
\end{eqnarray*}
where $\delta_{s,l} = 1$ if $s=l$, and zero otherwise.
Using the previous bounds, and the fact that $\operatorname{sc}_{0}(m,n) =
\operatorname{sc}_{0}(n,m)$, we obtain\vspace*{-1pt}
\begin{eqnarray*}
&&\bigl| \cov \bigl( P_{{2\pi s}/{T}}^{(T)}(m,n), P_{{2\pi
l}/{T}}^{(T)}(m,n)
\bigr) \bigr|
\\[-2pt]
&&\qquad \leq\frac{1}{4\pi^{2}} \bigl[ 4 T^{-2} \operatorname{sc}_{1}(m,m)
\operatorname{sc}_{1}(n,n) + 10 T^{-1} \operatorname{sc}_{0}(m,n,m,n)
\\[-2pt]
& &\hspace*{35pt}\qquad\quad{}+ 8 \delta_{s,l} \operatorname{sc}_{0}(m,m) \operatorname{sc}_{0}(n,n) + 8
\delta_{s+l, T} \operatorname{sc}_{0}(m,n)^{2} \bigr],
\end{eqnarray*}
and hence
%
\begin{eqnarray*}
T B_{T} \bigl| \var \bigl( \specdM{\omega}^{(T)}(m,n)
\bigr) \bigr|
&\leq& B_{T} \Biggl[ T^{-1} \sum
_{s = 1}^{T-1} W^{(T)}(\omega- 2\pi s/T)
\Biggr]^{2}
\\[-2pt]
&&\times\bigl[ 4 T^{-1} \operatorname{sc}_{1}(m,m)
\operatorname{sc}_{1}(n,n) + 10 \operatorname{sc}_{0}(m,n,m,n) \bigr]
\\[-2pt]
&&{}+8 \operatorname{sc}_{0}(m,m) \operatorname{sc}_{0}(n,n) B_{T}T^{-1}
\\[-2pt]
&&{}\times\sum_{s = 1}^{T-1} \bigl(W^{(T)}(
\omega- 2\pi s/T) \bigr)^{2}  + 8 \operatorname{sc}_{0}(m,n)^{2} B_{T} T^{-1}
\\[-2pt]
&&{}\times\sum_{s =
1}^{T-1} W^{(T)}(\omega- 2
\pi s/T)W^{(T)}(\omega+ 2\pi s/T).\vadjust{\goodbreak}
\end{eqnarray*}
Since at most $\frac{TB_{T}}{\pi} + 1$ of the summands are nonzero and
$\hnorm{W^{(T)}}_{\infty} \leq B_{T}^{-1} \hnorm{W}_{\infty}$
by Lemma~F.11, we obtain
$ [ T^{-1} \sum_{s = 1}^{T-1} W^{(T)}(\omega- 2\pi s/T)
]^{2} \leq\pi^{-2} \hnorm{W}_{\infty}^{2},$
and $B_{T}T^{-1} \sum_{s = 1}^{T-1} (W^{(T)}(\omega- 2\pi s/T))^{2}
\leq\pi^{-1} \hnorm{W}_{\infty}^{2},$
for large $T$. Similarly
$ \llvert  B_{T} T^{-1} \sum_{s = 1}^{T-1} W^{(T)}(\omega- 2\pi
s/T)W^{(T)}(\omega+ 2\pi s/T) \rrvert  \leq\pi^{-1} \hnorm{W}_{\infty}^{2}$
for large $T$. Since $B_{T} \rightarrow0,$ for $T$ large enough, we have
\begin{eqnarray*}
T B_{T} \bigl\llvert \var \bigl( \specdM{\omega}^{(T)}(m,n)
\bigr) \bigr\rrvert &\leq&\hnorm{W}_{\infty}^{2} \cdot \bigl[
\operatorname{sc}_{0}(m,n,m,n) + \operatorname{sc}_{1}(m,m)\operatorname{sc}_{1}(n,n)
\\
&&\hspace*{40pt}{} + 8 \operatorname{sc}_{0}(m,m)\operatorname{sc}_{0}(n,n) + 8 \operatorname{sc}_{0}(m,n)^{2}
\bigr].
\end{eqnarray*}
Now \eqref{eqspecdM-explicit-var} follows immediately by
setting $K = 8 \hnorm{W}_{\infty}^{2}.$

To prove \eqref{eqbound-sum-var-specdKE}, notice that, for large $T$,
inequality \eqref{eqspecdM-explicit-var} gives us
\begin{eqnarray*}
&&\sum_{m,n \geq1}T B_{T} \var \bigl( \specdM{
\omega}^{(T)}(m,n) \bigr) \\
&&\qquad\leq K \biggl[ \sum
_{m,n \geq1} \operatorname{sc}_{0}(m,n,m,n) + \biggl(\sum
_{m \geq
1}\operatorname{sc}_{1}(m,m) \biggr)^{2}
 + \biggl(\sum_{m \geq1} \operatorname{sc}_{0}(m,m)
\biggr)^{2} \\
&&\hspace*{223pt}\qquad{}+ \sum_{m,n \geq1}
\operatorname{sc}_{0}(m,n)^{2} \biggr].
\end{eqnarray*}
Notice that
$\cum{\xi_{t_{1}}(m), \xi_{t_{2}}(n), \xi_{t_{3}}(m), \xi_{0}(n)} = \sc
{\covO_{t_{1}, t_{2}, t_{3}} \vfi_{m} \otimes\vfi_{n}, \vfi_{m}
\otimes\vfi_{n}},$
hence
$ \sum_{m,n \geq1}\operatorname{sc}_0(m,n,m,n)
\leq \sum_{t_{1}, t_{2}, t_{3} \in\mathbb Z} \snorm{ \covO_{t_{1},
t_{2}, t_{3}} }_{1}.$
We also have\break $\cum{\xi_{t}(m),  \xi_{0}(n)} = \sc{\covO_{t} \vfi_{n},
\vfi_{m}},$
hence
$ \sum_{m \geq1}\operatorname{sc}_0(m,m) \leq\sum_{t \in\mathbb Z} \snorm{\covO_{t}}_{1}.$
Using the Cauchy--Schwarz inequality and Parseval's identity, we also obtain
$ \sum_{m,n \geq1}\operatorname{sc}_0(m,n)^{2}
\leq ( \sum_{t \in\mathbb Z} \snorm{\covO_{t}}_{1}  )^{2}.$
Similarly,
$\sum_{m,n \geq1} \operatorname{sc}_{1}(m,m) \operatorname{sc}_{1}(n,n) \leq  (\sum_{t \in
\mathbb Z} |t| \snorm{\covO_{t}}_{1}  )^{2}.$
Inequality \eqref{eqbound-sum-var-specdKE} is then established by
noticing that both $\sum_{t} \snorm{\covO_{t}}_{1}$ and $\sum_{t} |t|
\snorm{\covO_{t}}_{1}$ are bounded by
$\sum_{t \in\mathbb Z} (1 + |t|) \snorm{\covO_{t}}_{1},$
and setting $C = 3K.$

We can now put \eqref{eqspecdM-explicit-var} and \eqref
{eqbound-sum-var-specdKE} to use in order to establish the main
result. We first show that
$\sqrt{T B_{T}}  ( \specdK{\omega_{j}}^{(T)} - \ee \specdK{\omega
_{j}}^{(T)}  )$
is tight. Choose an orthonormal basis $\vfi_{n}$ of $ L^2( [0,1], \mathbb{R} )$. Notice that
\[
\mathbb{E} \bigl[
\bigl\langle \sqrt{T B_{T}} \bigl( \specdK{
\omega_{j}}^{(T)} - \ee \specdK {\omega_{j}}^{(T)}
\bigr), \vfi_{m} \otimes\vfi_{n} \bigr\rangle
^{2} \bigr] =
T B_{T} \var \bigl( \specdM{\omega}^{(T)}(m,n) \bigr).
\]
Since $(\vfi_{m} \otimes\vfi_{n})_{n,m \geq1}$ is an orthonormal
basis of $L^{2} ([0,1]^{2}, \mathbb C  )$, the tightness of
$\sqrt{T B_{T}}  ( \specdK{\omega_{j}}^{(T)} - \ee \specdK{\omega
_{j}}^{(T)}  )$ follows from \eqref{eqspecdM-explicit-var}, \eqref
{eqbound-sum-var-specdKE} and Lemma~\ref
{lmatech-criterion-unif-tightness}. Therefore the vector
$\sqrt{T B_{T}}  ( \specdK{\omega_{1}}^{(T)} - \ee \specdK{\omega
_{1}}^{(T)}, \ldots, \specdK{\omega_{J}}^{(T)} - \ee \specdK{\omega
_{J}}^{(T)}  )$
is also tight in $ (L^{2} ([0,1]^{2},\break \mathbb C  )
)^{J}.$ Applying \citet{Brill01}, Theorem 7.4.4, to the finite-dimensional
distributions of this vector completes the proof.
\end{pf}

Note here that condition (i) for $k=2$ is $\sum_{t \in\mathbb Z}
\hnorm{\covK_t}_2 < \infty,$ which guarantees that the spectral density
operator is continuous in $\omega$ with respect to the Hilbert--Schmidt
norm. If in addition we\vadjust{\goodbreak} want it to be continuous in $\tau, \sigma$ we
need to assume the stronger conditions $\sum_{t \in\mathbb Z} \hnorm
{\covK_t}_\infty< \infty$, and that each $\covK_t$ is continuous.

When $\omega=0$, the operator $2\pi\mathscr{F}_{\omega}$ reduces to the
long-run covariance operator $\sum_{t\in\mathbb{Z}}\mathscr{R}_t$, the
limiting covariance operator of the empirical mean. Correspondingly,
$2\pi\mathscr{F}^{(T)}_{0}$ is an estimator of the long-run covariance
operator that is consistent in mean square for the long-run covariance,
under no structural modelling assumptions. A similar estimator was also
considered in \citet{Horvath11}, who derived weak consistency under
$L^p$-$m$-approximability weak dependence conditions. \citet
{Hormann10} studied this problem by projecting onto a
finite-dimensional subspace. However, neither of these papers considers
functional central limit theorems for the estimator of the long-run
covariance operator; taking $\omega=0$, in Theorem \ref
{thmCLTspecdensity}, we obtain such a result:
%
\begin{cor}\label{corlimitlongrun}
Under the conditions of Theorem~\ref{thmCLTspecdensity}, we have
\[
\sqrt{B_T T} \bigl(2\pi\mathscr{F}_0^{(T)}-2
\pi\mathbb{E}\mathscr {F}_0^{(T)} \bigr)\stackrel{d} {
\longrightarrow}\mathcal{N} \bigl(0,(2\pi )^{3/2}\|W\|^2_2
\mathfrak{C} \bigr),
\]
where $\mathfrak{C}$ is the integral operator on $L^2([0,1]^2,\mathbb
{R})$ with kernel
\[
c(\tau_1,\sigma_1,\tau_2,
\sigma_2)= \bigl\{ \specdK{0}(\tau_{1}, \tau
_{2}) \specdK{0}(\sigma_{1}, \sigma_{2}) +
\specdK{0}(\tau_{1}, \sigma_{2}) \specdK{0}(
\sigma_{1}, \tau_{2}) \bigr\}.
\]
\end{cor}
We remark that the limiting Gaussian random operator is
purely real.

\section{Weak dependence, tightness and projections}\label{examples}

Our results on the asymptotic Gaussian representations of the discrete
Fourier transform and the spectral density estimator (Theorems \ref
{thmudf-asymptotics} and~\ref{thmCLTspecdensity}) effectively rest
upon two sets of weak dependence conditions: (1) the summability of the
nuclear norms of the autocovariance operators (at various rates), and
(2) the summability of the cumulant kernels of all orders (at various
rates). The roles of these two sets of weak dependence conditions are
distinct. The first is required in order to establish tightness of the
sequence of discrete Fourier transforms and spectral density estimators
of the underlying process. Tightness allows one to then apply the Cram\'
er--Wold device, and to determine the asymptotic distribution by
considering finite-dimensional projections; see, for example, \citet
{Ledo91}. The role of the second set of weak dependence conditions,
then, is precisely to allow the determination of the asymptotic law of
the projections, thus identifying the stipulated limiting distribution
via tightness.

Therefore, in principle, one can replace the second set of weak
dependence conditions with a set of conditions that allow for the
discrete Fourier transforms and spectral density estimators of the \emph
{vector} time series of the projections to be asymptotically Gaussian,
jointly in any finite number of frequencies. Our approach was to
generalise the cumulant multivariate conditions of \citet{Brill01},\vadjust{\goodbreak}
which do not require structural assumptions further to stationarity.
Alternatively, one may pursue generalizations of multivariate
conditions involving $\alpha$-mixing and summable cumulants of order
$2,4,$ and $8$ as in \citet{Hannan1970}, Chapter IV, Paragraph 4 and
\citeauthor{Rosenblatt1984} (\citeyear{Rosenblatt1984,Rosenblatt1985}), though $\alpha$-mixing can
also be a strong condition. Adding more structure, for example, in the
context of linear processes, one can focus on extending weaker
conditions requiring finite fourth moments and summable coefficients
[\citet{Hannan1970,Anderson1994}].

For the case of nonlinear moving-average representations of the form
$\xi_t = G(\vep_t, \vep_{t-1}, \ldots),$
where $G$ is a measurable function, and $\{\vep_j\}$ are i.i.d. random
variables, several results exist; however, none of them are (yet)
established for vector time series. For instance \citet{shao-wu2007}
show that if the second moment of $\xi_t$ is finite and
\[
\sum_{k = 0}^\infty\sqrt{ \ee \bigl\llvert
\mathbb{E}[\xi_k - \xi_{k+1} | \mathcal {F}_0]
\bigr\rrvert ^2} < \infty,
\]
where $\mathcal{F}_0$ is the sigma-algebra generated by $\{\vep_0, \vep
_{-1}, \ldots\}$, then the discrete Fourier transforms of $\xi_t$ are
asymptotically Gaussian, jointly for a finite number of frequencies.
Furthermore, \citet{shao-wu2007} establish the asymptotic normality of
the spectral density estimator at distinct frequencies under the moment
condition $\ee{|\xi_t|^{4+\delta}} < \infty,$ and the following
coupling condition: there exist $\alpha> 0$, $C > 0$ and $\rho\in
(0,1)$ such that
%
\begin{equation}
\label{eqcoupling-shao-wu} \ee\bigl| \xi_t - \xi_t'
|^\alpha< C \rho^t\qquad \forall t = 0,1,\ldots,
\end{equation}
where $\xi_t' = G(\vep_t, \ldots, \vep_1, \vep_0', \vep_{-1}', \ldots)$
and $(\vep_k')_{k \in\mathbb Z}$ is an i.i.d. copy of $(\vep_k)_{k \in
\mathbb Z}$.
Notice that \eqref{eqcoupling-shao-wu} is related to (in fact stronger
than) the $L^p$-$m$-approximability condition of \citet{Hormann10}.
Under the weaker conditions $\ee|\xi_t|^4 < \infty,$ and
\[
\sum_{t=0}^\infty \bigl( \ee|
\xi_t - \check\xi_t |^4 \bigr)^{1/4}
< \infty,
\]
where $\check\xi_t = G(\ldots, \vep_{-1}, \vep_{0}', \vep_{1}, \ldots,
\vep_{t})$ and $\vep_0'$ is an i.i.d. copy of $\vep_0,$ \citet
{liu-wu2010} establish that the spectral density estimator at a \emph
{fixed} frequency is asymptotically Gaussian. The idea behind these
coupling conditions is to approximate the series $\xi_t$ by
$m$-dependent series, for which derivation of asymptotic results is
easier. We also mention that, under milder conditions, \citet
{peligrad-wu2010} establish that for almost all $\omega\in(0,2\pi)$,
the discrete Fourier transform at $\omega$ is marginally asymptotically normal.

The weak dependence conditions pursued in this paper have the advantage
of not requiring additional structure, at the price of being relatively
strong if additional structure could be assumed. For example, if a
process is linear, the cumulant conditions will be satisfied provided
all moments exist and the coefficient operators are summable in an
appropriate sense, as shown in the proposition below. As mentioned\vadjust{\goodbreak}
above, we conjecture that four moments and summability of the
coefficients would suffice in the linear case; however, a more thorough
study of weak dependence conditions for the linear case is outside the
scope of the present paper.

\begin{prop}\label{propcumulant-decay-linear-process}
Let $X_{t} = \sum_{s \in\mathbb Z} A_{s}\vep_{t-s}$ be a linear
process with\break $\ee\hnorm{\vep_{0}}_{2}^{p} < \infty$ for all $p \geq
1$, and
$\sum_{s \in\mathbb Z} (1 + |s|^{l}) \hnorm{a_{s}}_{2} < \infty$ for
some positive integer $l$,
where $a_s$ is the kernel of $A_s$.
Then for all fixed $k=1,2, \ldots,$ $X_{t}$ satisfies {$\mathrm{C}(l,k)$},
%
\begin{eqnarray*}
\label{eqpropcumulant-decay-linear-process-0} \sum_{t_{1}, \ldots, t_{k-1} \in\mathbb Z}
\bigl(1 + |t_{j}|^{l} \bigr) \bigl\Vert \operatorname{cum}
(X_{t_{1}}, \ldots, X_{t_{k-1}}, X_{0} )\bigr \Vert_{2}
< \infty,\qquad
\forall j=1, \ldots, k-1.
\end{eqnarray*}
Furthermore,
\begin{eqnarray*}
\label{eqpropcumulant-decay-linear-process-1} \sum_{t \in\mathbb Z}
\bigl(1 + |t|^{l} \bigr) \vert\!\vert\!\vert \covO_{t} \vert\!\vert\!\vert_{1} &<& \infty,
\\
\label{eqpropcumulant-decay-linear-process-2} \sum_{t_{1}, t_{2}, t_{3}}
\bigl\vert\!\bigl\vert\!\bigl\vert \operatorname{cum} (X_{t_{1}}, X_{t_{2}},
X_{t_{3}}, X_{0} )\bigr \vert\!\bigr\vert\!\bigr\vert _{1} &< &\infty,
\end{eqnarray*}
where we view $\cum{X_{t_{1}}, X_{t_{2}}, X_{t_{3}}, X_{0}}$ as an
operator on $L^{2}([0,1]^{2}, \mathbb R)$; see Section~\ref
{secspectraldensity}.
\end{prop}

\begin{pf} See Proposition~E.1
in the supplementary material [\citet{supp}].
\end{pf}


\section{The effect of discrete observation}\label{discrete}

In practice, functional data are often observed on a discrete grid,
subject to measurement error, and smoothing is employed to make the
transition into the realm of smooth functions. This section considers
the stability of the consistency of our estimator of the spectral
density operator with respect to discrete observation of the underlying
stationary functional process. Since our earlier results do not a
priori require any smoothness of the functional data, except perhaps
smoothness that is imposed by our weak dependence conditions, we
consider a ``minimal'' scenario where the curves are only assumed to be
continuous in mean square. Under this weak assumption, we formalise the
asymptotic discrete observation framework via observation on an
increasingly dense grid subject to measurement error of variance
decreasing at a certain rate [e.g., \citet{hall2006b}]. {In principle,
one may drop the assumption that the noise variance decreases at a
certain rate at the expense of smoothness assumptions on the curves
that would suffice for smoothers constructed via the noisy sampled
curves to converge to the true curves, at a corresponding mean squared
error rate.}

Let $\Gamma$ be the grid $0=\tau_1 < \tau_2 < \cdots< \tau_M < \tau
_{M+1}=1$ on $[0,1]$, with $M=M(T)$ being a function of $T$ such that
$M(T) \rightarrow\infty$ as $T \rightarrow\infty,$ and
\[
|\Gamma| = \sup_{j=1, \ldots, M+1} \tau_{j} -
\tau_{j-1} \rightarrow0,\qquad M \rightarrow\infty.\vadjust{\goodbreak}
\]
Assume we observe the curves $X_t$ on this grid (except possibly at
$\tau_{M+1}$), additively corrupted by measurement error, represented
by independent and identically distributed random variables $\{
\varepsilon_{tj}\}$ (and independent of the $X_t$'s),
\[
y_{tj} = X_t(\tau_j) + \varepsilon_{tj},
\]
with $\ee\varepsilon_{tj} = 0$\vspace*{-1pt} and $\sqrt{\var(\varepsilon_{tj})} =
\sigma(M)$. Our goal is to show that our estimator of $\specdO{\omega
}^{(T)}$, when constructed on the basis of the $y_{tj}$'s, retains its
consistency for the true spectral density operator. To construct our
estimator on the basis of discrete observations, we use the following
(naive) proxy of the true~$X_t$,
\[
\lsub{\vep, s} {X_t(\tau)} = y_{tj}\qquad \mbox{if }
\tau_j \leq\tau < \tau_{j+1},
\]
and define the step-wise version of $X_t$,
\[
_{s}{X_t(\tau)} = X_t(\tau_j)\qquad
\mbox{if } \tau_j \leq\tau< \tau_{j+1}.
\]
Just as the spectral density kernel estimator $\specdK{\omega}^{(T)}$
is a functional of the $X_t$'s, we can define $\lsub{\vep, s}{\specdK
{\omega}^{(T)}}$ and $\lsub{s}{\specdK{\omega}^{(T)}}$, as the
corresponding functionals of the $\lsub{\vep,s}{X_t}$'s, $\lsub
{s}{X_t}$, respectively. The same can also be done for $\specdK{\omega
}$, $\specdO{\omega}^{(T)}$, $p_\omega^{(T)}$, $\udf{X}{\omega}^{(T)}$.
We then have the following stability result.

\begin{thmm} \label{thmdiscrete-to-continuous}
Under $\mathrm{C}(1,2)$, if $\ee\vep_{tj}^4 < \infty$, $\sigma^2(M) =
o(B_T),$ $B_T = o(1)$ such that $TB_T \rightarrow\infty,$
and if
%
\begin{equation}
\label{eqcov-infty-sum} \mbox{each } \covK_{t} \mbox{ is continuous,
and } \sum_{t} \Vert\covK_{t}\Vert_{\infty}
< \infty
\end{equation}
holds, then
\[
\int_{-\pi}^\pi \ee \bigl\vert\!\bigl\vert\!\bigl\vert \lsub{ \vep,
s} {\specdO{\omega}^{(T)} } - \specdO{\omega}^{(T)}\bigr \vert\!\bigr\vert\!\bigr\vert _2^2\ud\omega\rightarrow0, \qquad T \rightarrow\infty.
\]
Moreover,
we also have pointwise mean square convergence for a fixed $\omega$,
\[
\ee \bigl\vert\!\bigl\vert\!\bigl\vert \lsub{\vep, s} {\specdO{\omega}^{(T)} - \specdO{
\omega }^{(T)}} \bigr\vert\!\bigr\vert\!\bigr\vert ^{2}_{2}
\rightarrow0,\qquad T \rightarrow\infty
\]
under the same conditions if $0 < |\omega| < \pi,$ and under the
stronger condition $T B_{T}^{2} \rightarrow\infty$ if $\omega=0, \pm
\pi.$
\end{thmm}

\begin{pf*}{Proof of Theorem~\ref{thmdiscrete-to-continuous}}
First, we use the triangle inequality,
%
\begin{eqnarray}
\int_{-\pi}^\pi \ee \bigl\vert\!\bigl\vert\!\bigl\vert \lsub{ \vep,
s} {\specdO{\omega}^{(T)} } - \specdO{\omega}^{(T)}\bigr \vert\!\bigr\vert\!\bigr\vert _2^2 \,d\omega&=& \int_{-\pi}^\pi
\iint\ee \bigl\llvert \lsub{\vep, s} {\specdK{\omega}^{(T)} } - \specdK{
\omega}^{(T)} \bigr\rrvert ^2 \,d\omega
\nonumber
\\
&\leq& 2 \int_{-\pi}^\pi\iint\ee \bigl\llvert \lsub{
\vep, s} {\specdK {\omega}^{(T)} } - \lsub{s} {\specdK{
\omega}^{(T)} } \bigr\rrvert ^2 \,d\omega
\label{eqthm-discrete-to-cont-imse-eq1}
\\
& &{}+ 2 \int_{-\pi}^\pi\iint\ee \bigl\llvert \lsub
{s} {\specdK{\omega}^{(T)} } - \specdK{\omega}^{(T)} \bigr
\rrvert ^2 \,d\omega. \label{eqthm-discrete-to-cont-imse-eq2}
\end{eqnarray}
The inner integrals are on $[0,1]^2$ with respect to $d\tau \,d\sigma$.
First, we deal with the first summand,
\begin{eqnarray*}
\bigl\llvert \lsub{\vep, s} {\specdK{\omega}^{(T)} } - \lsub{s} {
\specdK{\omega }^{(T)} } \bigr\rrvert ^2 &=& 2\pi
T^{-2} \Biggl\llvert \sum_{l=0}^{T-1}
W^{(T)}(\omega- 2\pi l / T) \bigl( \lsub{\vep, s} {p_{2\pi l/T}^{(T)}}
- \lsub{s} {p_{2\pi l/T}^{(T)}} \bigr) \Biggr\rrvert ^2
\\
&\leq & O \bigl(T^{-1} \bigr) \sum_{l=0}^{T-1}
\bigl[W^{(T)}(\omega- 2\pi l / T) \bigr]^2 \bigl\llvert
\lsub{\vep, s} {p_{2\pi l/T}^{(T)}} - \lsub{s} {p_{2\pi
l/T}^{(T)}}
\bigr\rrvert ^2,
\end{eqnarray*}
where we have used Jensen's inequality.
We claim that, if $\tau_j \leq\tau< \tau_{j+1}$ and $\tau_k \leq
\sigma< \tau_{k+1}$,
\begin{eqnarray*}
\bigl\llvert \lsub{\vep, s} {p_\omega^{(T)}(\tau, \sigma)} -
\lsub{s} {p_\omega ^{(T)}(\tau, \sigma)} \bigr\rrvert
^2 &\leq&3 \bigl| \lsub{s} {\udf{X} {\omega }^{(T)}(\tau)}
\bigr|^2 \bigl|\udf{\vep} {-\omega}^{(T)}(k)\bigr|^2
\\
&&{}+ 3 \bigl| \udf{\vep} {\omega}^{(T)}(j) \udf{\vep } {-\omega}^{(T)}(k)\bigr|^2
+ 3 \bigl| \udf{\vep} {\omega}^{(T)}(j) \bigr|^2 \bigl| \lsub {s} {\udf{X}
{-\omega}^{(T)}(\sigma)} \bigr|^2,
\end{eqnarray*}
where $\udf{\vep}{\omega}^{(T)}(j) = (2\pi T)^{-1/2} \sum_{l=0}^{T-1}
e^{-\ii\omega t} \vep_{tj}.$ To see this, we note that
\begin{eqnarray*}
\lsub{\vep, s} {p_\omega^{(T)}(\tau, \sigma)} - \lsub{s}
{p_\omega ^{(T)}(\tau, \sigma)} &=& \lsub{\vep, s} {\udf{X} {
\omega}^{(T)}(\tau)} \cdot\lsub{\vep, s} {\udf{X} {-\omega}^{(T)}(
\sigma)} - \lsub{s} {\udf {X} {\omega}^{(T)}(\tau)} \cdot\lsub{s} {
\udf{X} {-\omega}^{(T)}(\sigma)}
\\
&=& \bigl( \lsub{\vep, s} {\udf{X} {\omega}^{(T)}(\tau)} - \lsub {s} {
\udf{X} {\omega}^{(T)}(\tau)} \bigr) \cdot\lsub{\vep, s} {\udf {X} {-
\omega}^{(T)}(\sigma)}
\\
&&{}+ \lsub{s} {\udf{X} {\omega}^{(T)}(\tau)} \cdot \bigl( \lsub{\vep, s}
{\udf{X} {-\omega}^{(T)}(\sigma)} - \lsub{s} {\udf {X} {-
\omega}^{(T)}(\sigma)} \bigr)
\\
&=& \lsub{s} {\udf{X} {\omega}^{(T)}(\tau)} \udf{\vep} {-\omega
}^{(T)}(k) + \udf{\vep} {\omega}^{(T)}(j) \udf{\vep} {-
\omega}^{(T)}(k)
\\
&&{} + \udf{\vep} {\omega}^{(T)}(j) \lsub{s} {\udf {X} {-
\omega}^{(T)}(\sigma)},
\end{eqnarray*}
since $\lsub{\vep, s}{\udf{X}{\omega}^{(T)}(\tau)} = \lsub{s}{\udf
{X}{\omega}^{(T)}(\tau)} + \udf{\vep}{\omega}^{(T)}(j)$, and similarly
if we replace $\sigma$ by $\tau$ and $j$ by $k$. Our claim thus follows
from Jensen's inequality.

In order to bound the expectation of $\llvert  \lsub{\vep, s}{p_\omega
^{(T)}(\tau, \sigma)} - \lsub{s}{p_\omega^{(T)}(\tau, \sigma)} \rrvert ^2$, we will first compute the expectation, conditional on the $\sigma
$-algebra generated by the $X_t$'s, which we will denote by $\ee_X$,
and then use the tower property. As an intermediate step, we claim that
$\ee_X | \udf{\vep}{\omega}^{(T)}(j) |^2 = O(\sigma^2(M)),$
\[
\ee_X \bigl| \udf{\vep} {\omega}^{(T)}(j) \udf{\vep} {-
\omega}^{(T)}(k) \bigr|^2 = \cases{ O \bigl( \sigma^4(M)
\bigr), & \quad$\mbox{if } j\neq k,$ \vspace*{2pt}
\cr
O \bigl(\sigma^4(M)
\bigr) + O \bigl(T^{-1} \bigr), &\quad $\mbox{if } j=k,$}
\]
uniformly in $j,k$ (notice that all $\ee_X$ can be replaced by $\ee$
since the $\vep_{tj}$'s are independent of the $X_t$'s). To establish
this, notice that $| \udf{\vep}{\omega}^{(T)}(j) |^2 = \udf{\vep}{\omega
}^{(T)}(j) \udf{\vep}{-\omega}^{(T)}(j),$ hence
$\ee_X | \udf{\vep}{\omega}^{(T)}(j) |^2 = (2\pi T)^{-1} \sum_{t,s=0}^{T-1} e^{-\ii\omega(t-s)} \eee{\vep_{tj} \vep_{sj}}.$
The summand is equal to $\sigma^2(M)$ if $t=s$, and zero otherwise (by
independence of the $\vep$'s), hence the first statement follows
directly. The case $j \neq k$ follows from the first statement, once
the independence of $\udf{\vep}{\omega}^{(T)}(j)$ and $\udf{\vep
}{-\omega}^{(T)}(k)$ has been noticed.
We can now turn to the case $j=k$. First notice that
\[
\ee_X \bigl| \udf{\vep} {\omega}^{(T)}(j) \udf{\vep} {-
\omega}^{(T)}(j)\bigr |^2 = (2\pi T)^{-2} \sum
_{t_1, t_2, t_3, t_4 = 0}^{T-1} e^{-\ii\omega[(t_1 - t_2) + (t_3 - t_4)]} \ee_X (
\vep_{t_1} \vep _{t_2} \vep_{t_3}
\vep_{t_4}),
\]
where we have written $\vep_{t}$ instead of $\vep_{t j}$ for tidiness.
The expectation of the product of the $\vep$'s is equal to zero if at
least one of the $t_l$'s is different from the all the other ones (by
independence). So we may assume that each $\vep_{t_l}$ appears at least
twice. There can be therefore $2-r$ distinct terms in $\vep_{t_1} \vep
_{t_2} \vep_{t_3} \vep_{t_4}$, where $r=0$ or $1$. If $r=0$,
$\ee_X (\vep_{t_1} \vep_{t_2} \vep_{t_3} \vep_{t_4}) = \sigma^4(M),$
and if $r=1$,
$\ee_X (\vep_{t_1} \vep_{t_2} \vep_{t_3} \vep_{t_4}) = \ee_X \vep^4 =
\ee\vep^4.$
Thus $\ee_X | \udf{\vep}{\omega}^{(T)}(j) \udf{\vep}{-\omega}^{(T)}(j)
|^2 = (2\pi T)^{-2} [ N_0 \sigma^4(M) + N_1 \ee\vep^4  ],$
where $N_r$ is the number of ways we can assign integers $t_1, \ldots,
t_4$ in $\{0, \ldots, T-1\}$ such that each $t_l$ appears at least
twice and exactly $2-r$ distinct integers appear. Simple combinatorics yield
$N_0 = \bigl({{4}\atop{2}} \bigr)T(T-1) = 6T(T-1),$
and $N_1 = T,$ and so the case case $j \neq k$ follows directly since
$\ee\vep^4 < \infty.$

We can now bound $\ee_X \llvert  \lsub{\vep, s}{p_{2\pi l/T}^{(T)}} -
\lsub{s}{p_{2\pi l/T}^{(T)}} \rrvert ^2$,
\begin{eqnarray*}
\ee_X \bigl\llvert \lsub{\vep, s} {p_{2\pi l/T}^{(T)}}
- \lsub{s} {p_{2\pi
l/T}^{(T)}} \bigr\rrvert ^2 &\leq&3
\bigl| \lsub{s} {\udf{X} {2\pi l/T}^{(T)}(\tau)} \bigr|^2
\ee_X \bigl|\udf{\vep} {-2\pi l/T}^{(T)}(k)\bigr|^2
\\
&&{} + 3 \ee_X \bigl| \udf{\vep} {2\pi l/T}^{(T)}(j) \udf { \vep}
{-2\pi l/T}^{(T)}(k)\bigr|^2
\\
&&{} + 3 \bigl| \lsub{s} {\udf{X} {-2\pi l/T}^{(T)}(\sigma )} \bigr|^2
\ee_X \bigl| \udf{\vep} {2\pi l/T}^{(T)}(j) \bigr|^2
\\
&\leq& O \bigl(\sigma^2(M) \bigr) \bigl[ \bigl| \lsub{s} {\udf{X} {2\pi
l/T}^{(T)}(\tau )}\bigr |^2 + \bigl| \lsub{s} {\udf{X} {-2\pi
l/T}^{(T)}(\sigma)} \bigr|^2 \bigr]
\\
&&{} + O \bigl( \sigma^4(M) \bigr) + O \bigl(T^{-1} \bigr).
\end{eqnarray*}
Since
$ | {\udf{X}{\omega}^{(T)}(\tau)} |^2 = | p_\omega^{(T)}(\tau, \tau)|,$
Proposition~\ref{propapprox-expectation-periodogram}, Remark~\ref
{rmksup-norm-cumulant-conditions} and \eqref{eqcov-infty-sum} yield
that $\int\ee| \lsub{s}{\udf{X}{2\pi l/T}^{(T)}(\tau)} |^2 \,d\tau=
O(1)$. Using the tower property, we obtain
\[
\iint\ee \bigl\llvert \lsub{\vep, s} {p_{2\pi l/T}^{(T)}} -
\lsub{s} {p_{2\pi
l/T}^{(T)}} \bigr\rrvert ^2 \leq O
\bigl(\sigma^2(M) \bigr) + O \bigl(T^{-1} \bigr),
\]
uniformly in $l=1, \ldots, T-1$ under the assumptions of this theorem. Thus
\begin{eqnarray*}
&&\iint\ee \bigl\llvert \lsub{\vep, s} {\specdK{\omega}^{(T)} } - \lsub
{s} {\specdK{\omega}^{(T)} } \bigr\rrvert ^2 \\
&&\qquad\leq O
\bigl(T^{-1} \bigr) \sum_{l=0}^{T-1}
\bigl[W^{(T)}(\omega- 2\pi l / T) \bigr]^2 \cdot\iint\ee
\bigl\llvert \lsub{\vep, s} {p_{2\pi l/T}^{(T)}} - \lsub{s}
{p_{2\pi l/T}^{(T)}} \bigr\rrvert ^2
\\
&&\qquad= O \bigl(B_T^{-1}\sigma^2(M) \bigr) +
O(B_T T)^{-1},
\end{eqnarray*}
uniformly in $ \omega.$ Hence we obtain the bound on the expectation of
first summand~\eqref{eqthm-discrete-to-cont-imse-eq1},
\[
\int_{-\pi}^\pi\iint\ee \bigl\llvert \lsub{\vep,
s} {p_{2\pi l/T}^{(T)}} - \lsub{s} {p_{2\pi l/T}^{(T)}}
\bigr\rrvert ^2 \,d\omega= O \bigl(B_T^{-1}
\sigma^2(M) \bigr) + O(B_T T)^{-1},
\]
under the assumptions of the theorem.

We now turn to the second summand \eqref
{eqthm-discrete-to-cont-imse-eq2}. First notice that
\[
\int_{-\pi}^\pi\iint\ee \bigl\llvert \lsub{s} {
\specdK{\omega}^{(T)} } - \specdK{\omega}^{(T)} \bigr\rrvert
^2 \,d\omega= 2 \int_{0}^\pi\iint\ee
\bigl\llvert \lsub{s} {\specdK{\omega}^{(T)} } - \specdK{
\omega}^{(T)} \bigr\rrvert ^2 \,d\omega,
\]
since $ \lsub{s}{\specdK{-\omega}^{(T)}} = \overline{ \lsub{s}{\specdK
{\omega}^{(T)}} }$ and ${\specdK{-\omega}^{(T)}} = \overline{ {\specdK
{\omega}^{(T)}} }$.
Using the decomposition
\[
\ee \bigl\llvert \lsub{s} {\specdK{\omega}^{(T)} } - \specdK{
\omega}^{(T)} \bigr\rrvert ^2 = \cov \bigl[ \lsub{s} {
\specdK{\omega}^{(T)} } - \specdK {\omega}^{(T)}, \lsub{s} {
\specdK{\omega}^{(T)} } - \specdK{\omega }^{(T)} \bigr] + \bigl
\llvert \mathbb{E} \bigl[\lsub{s} {\specdK{\omega}^{(T)} } - \specdK{
\omega}^{(T)} \bigr] \bigr\rrvert ^2,
\]
the covariance term can be written as sums and differences of four
terms of the form $\cov( \specdK{\omega}^{(T)}(\sigma_1, \sigma_2),
\specdK{\omega}^{(T)}(\sigma_3, \sigma_4))$, for some $\sigma_l$'s. The
important thing here is that each of these terms can be bounded in
$L^{2}$---independently of the $\sigma_l$'s---using Corollary~\ref
{corcov-unif-bound} and Proposition~\ref{propcov-sharp-bound},
\begin{eqnarray*}
&&\cov \bigl[ \lsub{s} {\specdK{\omega}^{(T)} } - \specdK{
\omega}^{(T)}, \lsub{s} {\specdK{\omega}^{(T)} } - \specdK{
\omega}^{(T)} \bigr] \\
&&\qquad= \cases{ O \bigl(B_T^{-2}T^{-1}
\bigr) + O \bigl(T^{-1} \bigr), &\quad $\mbox{if } \omega\in[0,
B_T]\cup [\pi- B_T, \pi],$ \vspace*{2pt}
\cr
O
\bigl(B_T^{-1}T^{-1} \bigr), & \quad$\mbox{if } \omega
\in[B_T,\pi- B_T]$ }
\end{eqnarray*}
in $L^{2}$. Hence decomposing $\int_0^\pi= \int_0^{B_T} + \int_{B_T}^{\pi- B_T} + \int_{\pi- B_T}^{\pi}$, we obtain
\[
\int_{0}^\pi\iint\cov \bigl[ \lsub{s} {\specdK{
\omega}^{(T)} } - \specdK{\omega}^{(T)}, \lsub{s} {\specdK{
\omega}^{(T)} } - \specdK{\omega }^{(T)} \bigr] \,d\omega= O
\bigl(B_T^{-1} T^{-1} \bigr),
\]
if $B_T \rightarrow0$.

In order to bound $\llvert  \eee{ \lsub{s}{\specdK{\omega}^{(T)} } -
\specdK{\omega}^{(T)}} \rrvert ^2$, we use Proposition~\ref
{propexpectation-specdk-estimator} and Lemma~F.4 (with $p=1$),
\[
\iint \bigl\llvert \mathbb{E} \bigl[\lsub{s} {\specdK{\omega}^{(T)} } -
\specdK{\omega }^{(T)} \bigr] \bigr\rrvert ^2 \leq4 \iint
\llvert \lsub{s} {\specdK{\omega}} - \specdK{\omega} \rrvert ^2 + O
\bigl(B_T^2 \bigr) + O \bigl(T^{-2} \bigr) +
O(B_T T)^{-2},
\]
uniformly in $\omega.$ Thus
\begin{eqnarray*}
&&\int_{0}^\pi\iint \bigl\llvert \mathbb{E} \bigl[
\lsub{s} {\specdK{\omega}^{(T)} } - \specdK{\omega}^{(T)}
\bigr] \bigr\rrvert ^2 \,d\omega\\
&&\qquad\leq4 \int_{-\pi}^\pi
\iint \llvert \lsub{s} {\specdK{\omega}} - \specdK{\omega} \rrvert ^2
\,d\omega+ O \bigl(B_T^2 \bigr) + O
\bigl(T^{-2} \bigr) + O(B_T T)^{-2}.
\end{eqnarray*}
The quantity $\iint\llvert  \lsub{s}{\specdK{\omega}} - \specdK{\omega}
\rrvert ^2$ is in fact the the squared distance between $\lsub{s}{\specdK
{\omega}}$ and $\specdK{\omega}$ in the space $L^2([0,1]^2, \mathbb
C)$. Under \eqref{eqcov-infty-sum}, $\specdK{\omega}(\tau, \sigma)$ is
uniformly continuous in $\omega, \tau, \sigma$; since $\lsub{s}{\specdK
{\omega}}$ is a step-wise approximation of $\specdK{\omega}$, we obtain
\[
\sup_{\omega\in[-\pi, \pi]} \iint\llvert \lsub{s} {\specdK{\omega}} - \specdK{
\omega} \rrvert ^2 \rightarrow0,\qquad M \rightarrow\infty.
\]
Piecing these results together, we obtain
\[
\int_{0}^\pi\iint\ee \bigl\llvert \lsub{s} {
\specdK{\omega}^{(T)} } - \specdK{\omega}^{(T)} \bigr\rrvert
^2 \,d\omega= o(1) + O \bigl(B_T^{-1}
T^{-1} \bigr) + O \bigl(B_T^2 \bigr),
\]
and therefore
\[
\int_{-\pi}^\pi\iint\ee \bigl\llvert \lsub{\vep,
s} {\specdK{\omega}^{(T)} } - \specdK{\omega}^{(T)} \bigr
\rrvert ^2 \,d\omega= O \bigl(\sigma^2(M)B_T^{-1}
\bigr) + O \bigl(B_T^{-1} T^{-1} \bigr) + O
\bigl(B_T^2 \bigr) + o(1),
\]
where the $o(1)$ term comes from the $L^2$ distance between $\lsub
{s}{\specdK{\omega}}$ and $\specdK{\omega}$. Under our assumptions, the
right-hand side tends to zero as $T\rightarrow\infty$.

A careful examination of the proof reveals that the
pointwise statement of the theorem follows with a directly analogous argument.
\end{pf*}

\begin{rmk} \label{rmkdiscrete-to-cont-best-rate}
The use of Proposition~\ref{propcov-sharp-bound} was valid in this
context, but requires some attention. Indeed, it relies on Lemma~F.15 in the supplementary material [\citet{supp}], applied to
$g_{(\tau, \sigma)}(\alpha) = \lsub{s}{\specdK{\alpha}^{(T)}(\tau,
\sigma)}.$ Remark~F.16 in the supplementary material [\citet{supp}] tells us that the convergence of the convolution integral
depends on the {uniform} continuity parameter $\delta(\vep)$, which
here will depend on the size of the sampling grid $M=M(T)$; in other
words, $\delta(\vep) = \delta(\vep, M)$. But notice that since \eqref
{eqcov-infty-sum} holds,
\begin{eqnarray*}
\Vert \lsub{s} {\specdK{\omega_1}} - \lsub{s} {\specdK{
\omega_2}} \Vert_{2} &\leq&\sup_{0 \leq\tau, \sigma\leq1} \bigl
\llvert \lsub{s} {\specdK {\omega_1}(\tau, \sigma)} - \lsub{s} {
\specdK{\omega_2}(\tau, \sigma)} \bigr\rrvert
\\
&=& \sup_{\tau, \sigma= \tau_1, \ldots, \tau_M} \bigl\llvert {\specdK {\omega_1}(
\tau, \sigma)} - {\specdK{\omega_2}(\tau, \sigma)} \bigr\rrvert
\\
&\leq& \sup_{0 \leq\tau, \sigma\leq1} \bigl\llvert {\specdK{\omega
_1}(\tau, \sigma)} - {\specdK{\omega_2}(\tau, \sigma)}
\bigr\rrvert,
\end{eqnarray*}
hence we can choose a $\delta(\vep)$ that is independent of $M$, and
the application of Proposition~\ref{propcov-sharp-bound} is valid.
\end{rmk}


\section{Numerical simulations}\label{simulations}

In order to probe the finite sample performance of our estimators (in
terms of IMSE), we have performed numerical simulations on stationary
functional time series admitting a linear representation
\[
X_{t} = \sum_{s = 0}^{10}
A_{s} \vep_{t-s}.
\]
We have taken the collection of innovation functions $\{\vep_{t}\}$ to
be independent Wiener processes on $[0,1]$, which we have represented
using a truncated Karhunen--Lo\`eve expansion,
\[
\vep_{t}(\tau) = \sum_{k=1}^{1000}
\xi_{k,t} \sqrt{\lambda_{k}} e_{k}(\tau).
\]
Here $\lambda_{k} = 1/[(k - 1/2)^{2} \pi^{2}]$, $\xi_{k,t}$ are
independent standard Gaussian random variables and
$e_{k}(\tau) = \sqrt{2} \sin[(k - 1/2) \pi\tau]$
is orthonormal system in $L^{2}([0,1], \mathbb R)$ [\citet
{adler1990}]. We have constructed the operators $A_{s}$ so that their
image be contained within a $50$-dimensional subspace of $ L^2( [0,1], \mathbb{R} )$,
spanned by an orthonormal basis $\psi_{1}, \ldots, \psi_{50}$.
Representing $\vep_{t}$ in the $e_{k}$ basis, and $A_{s}$ in the $\psi
_m \otimes e_{k}$ basis, we obtain a matrix representation of the
process $X_{t}$ as
$\b{X}_{t} = \sum_{s=0}^{10} \b{A}_{s} \bvep_{t-s},$
where $\b{X}_{t}$ is a $50 \times1$ matrix, each $\b{A}_{s}$ is a $50
\times1000$ matrix, and each $\bvep_{t}$ is a $1000 \times1$ matrix.

We simulated a stretch of $X_{t}, t=0, \ldots, T-1$ for $T=2^{n}$, with
$n=7,8, \ldots,\break 15.$ Typical functional data sets would range between
$T=2^6$ and $T=2^8$ data points. We constructed the matrices $\b
{A}_{s}$, as random Gaussian matrices with independent entries, such
that elements in row $j$ where $N(0, j^{-2\alpha})$ distributed. When
$\alpha=0$, the projection of each $\vep_t$ onto the subspace spanned
by each $\psi_m, m=1, \ldots, 50$ has (roughly) a comparable magnitude.
A positive value of $\alpha$, for example, $\alpha=1$ means that the
projection of $\vep_t$ onto the subspace spanned by $\psi_j$ will have
smaller magnitude for larger $j$'s.

For comparison purposes, we also carried out analogous simulations,
but with $\lambda_{k} = 1$, that is, the variance of the innovations
$\vep_{t}$ being equal to one in each direction $e_{n}, n=1, \ldots,
1000.$ In the sequel, we will refer to these as the simulations with
``white noise innovations,'' and to the previous ones as ``Wiener
innovations.'' The white noise process is, of course, not a true white
noise process, but a projection of a white noise process. However, it
does represent a case of a ``rough'' innovation process, which we present
here as an extreme scenario.

For each $T$, we generated 200 simulation runs which we used to compute
the IMSE by approximating the integral
\[
2\int_{0}^{\pi} \ee \bigl\vert\!\bigl\vert\!\bigl\vert \specdO{
\omega} - \specdO{\omega }^{(T)} \bigr\vert\!\bigr\vert\!\bigr\vert _{2}^{2}
\,d \omega
\]
by a weighted sum over the finite grid $\Gamma= \{ \pi j / 10;
j=0, \ldots, 9 \}.$ We chose $B_{T} = T^{-1/5}$ [e.g., \citet{Grenander1966}, Paragraph
4.7, \citet{Brill01}, Paragraph 7.4] and $W(x)$ to be
the \emph{Epanechnikov kernel} [e.g., \citet{Wand1995}], $W(x) = \frac
{3}{4}(1-x^2)$ if $|x| < 1,$ and zero otherwise. The results are shown
in a log-log scale in Figure~\ref{figimse}, for $\alpha=2$. The slopes
of the least square lines passing through the medians of the simulation
results show that $\operatorname{IMSE}(\specdO{}^{(T)}) \propto T^{\beta}$, with
$\beta\approx-0.797$ for the white noise innovations, and $\beta
\approx-0.796$ for the Wiener innovations. According to Theorem~\ref
{thmimse-convergence}, the decay of the $\operatorname{IMSE}(\specdO{}^{(T)})$ is
bounded by
\[
C_{1} T^{-2/5} + C_{2} T^{-4/5} \approx
C_{1} T^{-0.4} \qquad\mbox{(if $T$ is large)}
\]
for some constants $C_{1}, C_{2}$.

\begin{figure}

\includegraphics{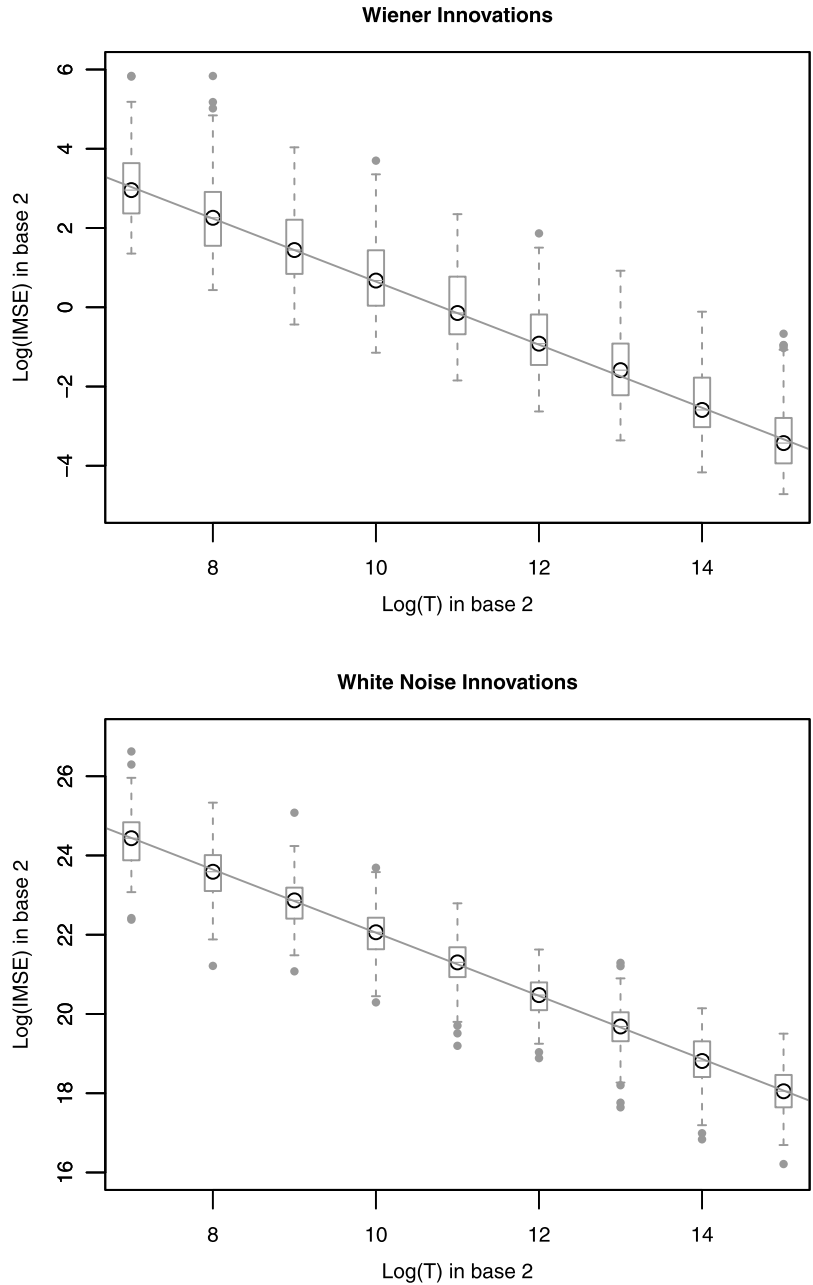}

\caption{The results of the simulated ISE in a log--log scale, {with
$\alpha= 2$.
The upper and lower plot correspond to the Wiener Innovations and the
White Noise Innovation setups, respectively}.
The dots correspond to the median of the results of the simulations,
and the lines are the least square lines of the medians.
The boxplots summarise the distribution of the ISE for the 200
simulation runs. Though the ranges of the $y$-axes are different,
the scales are the same, and the two least square lines are indeed
almost parallel.} \label{figimse}
\end{figure}

\begin{figure}

\includegraphics{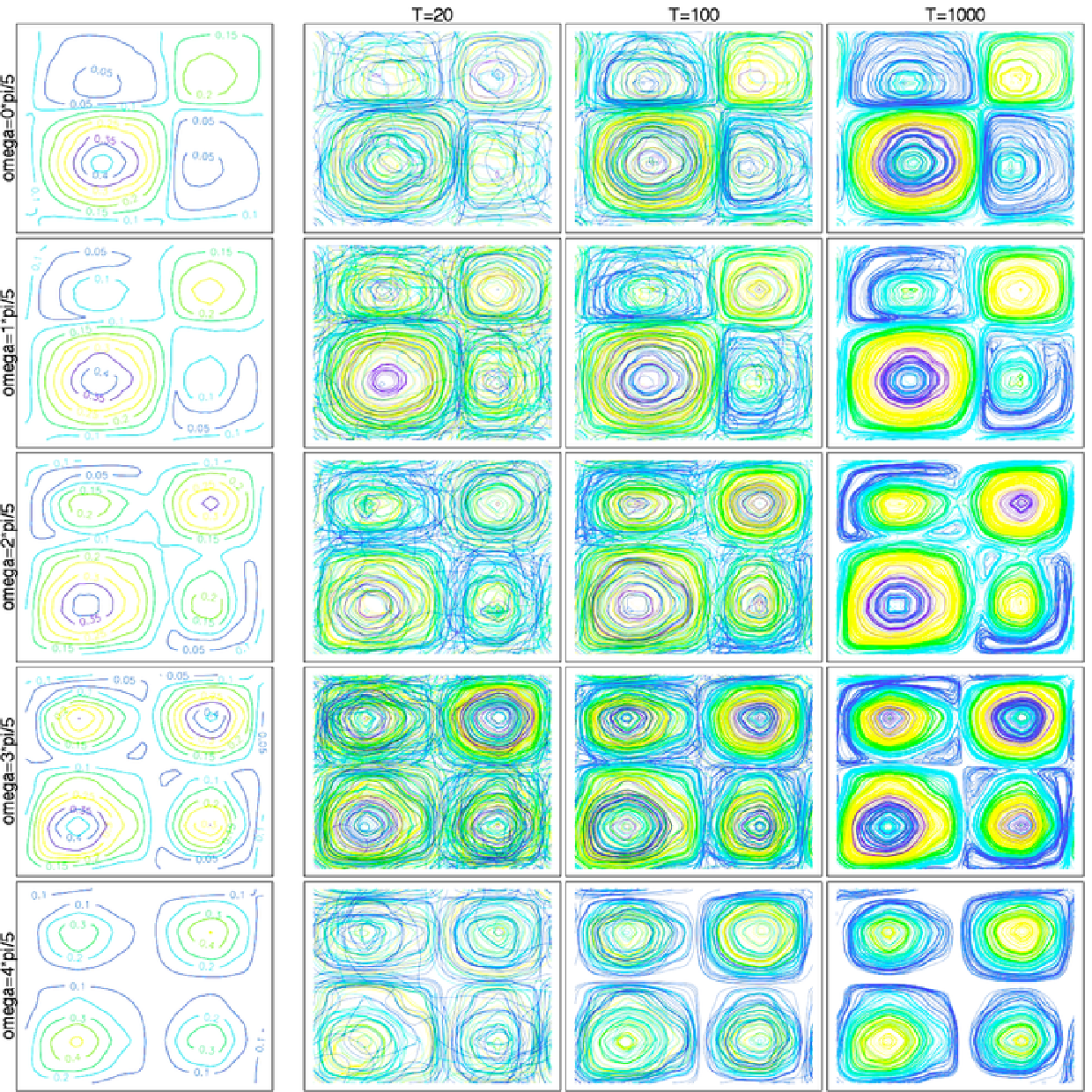}

\caption{Contour plots for the amplitude of the true and
estimated spectral density kernel when the innovation
process consists of Wiener processes. Each row corresponds
to a different frequency ($\omega=k\pi/5$, $k=0,1,\ldots,4$,
going from top to bottom). The first column contains the
contour plots of the true amplitudes of the kernel at each
corresponding frequency. The rest of the columns correspond
to the estimated contours for different sample sizes ($T=20, 100, 1000$
from left to right).
Twenty estimates, corresponding to twenty replications
of the process, have been superposed in order to provide
a visual illustration of the variability. The contours
plotted always correspond to the same level curves and
use the same colour-coding in each row.} \label{figcontourWiener}
\end{figure}

\begin{figure}

\includegraphics{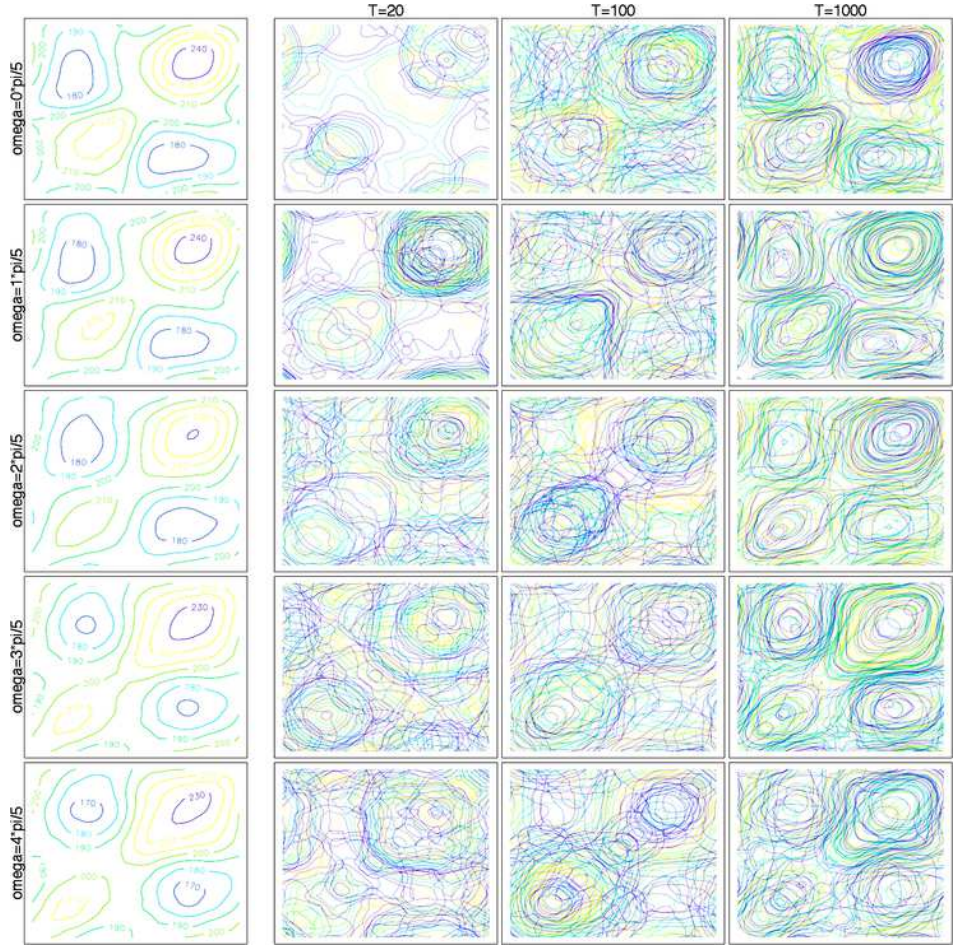}

\caption{Contour plots for the amplitude of the true and estimated
spectral density kernel when the innovation process consists of white
noise processes. Each row corresponds to a different frequency
($\omega=k\pi/5$, $k=0,1,\ldots,4$, going from top to bottom).
The first column contains the contour plots of the true
amplitudes of the kernel at each corresponding frequency.
The rest of the columns correspond to the estimated contours
for different sample sizes ($T=20, 100, 1000$ from left to right).
Twenty estimates, corresponding to twenty replications of the process,
have been superposed in order to provide a visual illustration of the
variability.
The contours plotted always correspond to the same level curves and use
the same colour-coding in each row.}\label{figcontourwhite}
\end{figure}

In order to gain a visual appreciation of the accuracy of the
estimators, we construct plots to compare the true and estimated
spectral density kernels in Figures~\ref{figcontourWiener} and \ref
{figcontourwhite}, for the Wiener and white noise cases, respectively.
For practical purposes, we set $\alpha=2$, as for the simulation of
the IMSEs. We simulated $X_t = A_0 \vep_t + A_1 \vep_{t-1}$, where $\vep
_t(\tau)$ lies on the subspace of $L^2([0,1], \mathbb R)$ spanned by
the basis $e_1, \ldots, e_{100}$, and the operators $A_0, A_1$ lie in
the subspace spanned by $(\psi_m \otimes e_k)_{m=1,\ldots, 51; k=1,
\ldots, 100}$.
Since the target parameter is a complex-valued function defined over a
two-dimensional rectangle, some information loss must be incurred when
representing it graphically. We chose to suppress the phase component
of the spectral density kernel, plotting only its amplitude, $|\specdK
\omega(\tau,\sigma)|$, for all $(\tau,\sigma)\in[0,1]^2$ and for
selected frequencies $\omega$ (the spectral density kernel is seen to
be smooth in $\omega$, so this does not entail a significant loss of
information). For various choices of sample size~$T$, we have
replicated the realisation of the process, and the corresponding kernel
density estimator for the particular frequency. Each time, we plotted
the contours in superposition, in order to be able to visually
appreciate the variability in the estimators: tangled contour lines
where no clear systematic pattern emerges signify a region of high
variability, whereas aligned contour lines that adhere to a
recognisable shape represent regions of low variability. As is
expected, the ``smoother'' the innovation process, the less variable the
results appear to be, and the variability decreases for larger values
of $T$.


\section{Background results and technical statements}\label{technical}

Statements and proofs of intermediate results in functional analysis
and probability in function space that are required in our earlier
formal derivations, can be found in the supplementary material [\citet{supp}]. This
supplement also collects some known results and facts for the reader's
ease. We include here a useful lemma that provides an easily verifiable
$L^2$ moment condition that is sufficient for tightness to hold true.
It collects arguments appearing in the proof of \citet{Bosq00}, Theorem~2.7, and its proof can also be found in the supplementary material [\citet{supp}].

\begin{lma}[(Criterion for tightness in Hilbert space)] \label
{lmatech-criterion-unif-tightness}
Let $H$ be a (real or complex) separable Hilbert space, and $X_{T} \dvtx \Omega\rightarrow H,   T=1,2,\ldots,$ be a sequence of random
variables. If for some complete orthonormal basis $\{e_{n}\}_{n\geq1}$
of $H$, we have $\ee\llvert  \sc{X_{T}, e_{n}} \rrvert ^{2} \leq a_{n},
n=1,2, \ldots,$
for all large $T$, and
$\sum_{n \geq1} a_{n} < \infty,$
then $\{X_{T}\}_{T\geq1}$ is tight.

\end{lma}

\section*{Acknowledgements}
Our thanks go the Editor, Associate Editor and three referees for
their careful reading and thoughtful comments.

\begin{supplement}[id=suppA]\label{supplement}
\sname{Online Supplement}
\stitle{``Fourier Analysis of Stationary Time Series in Function Space''}
\slink[doi]{10.1214/13-AOS1086SUPP} 
\sdatatype{.pdf}
\sfilename{aos1086\_supp.pdf}
\sdescription{The online supplement
contains the proofs that were omitted, and several additional technical
results used in this paper.}
\end{supplement}


\begin{thebibliography}{59}

\bibitem[\protect\citeauthoryear{Adler}{1990}]{adler1990}
%
\begin{bbook}[mr]
\bauthor{\bsnm{Adler},~\bfnm{Robert~J.}\binits{R.~J.}}
(\byear{1990}).
\btitle{An Introduction to Continuity, Extrema, and Related Topics for General
{G}aussian Processes}.
\bseries{Institute of Mathematical Statistics Lecture Notes---Monograph Series}
\bvolume{12}.
\bpublisher{IMS}, \blocation{Hayward, CA}.
\bid{mr={1088478}}
\bptok{imsref}%
\end{bbook}
%
\endbibitem

\bibitem[\protect\citeauthoryear{Anderson}{1994}]{Anderson1994}
%
\begin{bbook}[mr]
\bauthor{\bsnm{Anderson},~\bfnm{T.~W.}\binits{T.~W.}}
(\byear{1994}).
\btitle{The Statistical Analysis of Time Series}.
\bpublisher{Wiley}, \blocation{New York}.
\bptnote{check year}%
\bptok{imsref}%
\end{bbook}
%
\endbibitem

\bibitem[\protect\citeauthoryear{Antoniadis, Paparoditis and
Sapatinas}{2006}]{antoniadis2006}
%
\begin{barticle}[mr]
\bauthor{\bsnm{Antoniadis},~\bfnm{Anestis}\binits{A.}},
\bauthor{\bsnm{Paparoditis},~\bfnm{Efstathios}\binits{E.}} \AND
\bauthor{\bsnm{Sapatinas},~\bfnm{Theofanis}\binits{T.}}
(\byear{2006}).
\btitle{A functional wavelet-kernel approach for time series prediction}.
\bjournal{J. R. Stat. Soc. Ser. B Stat. Methodol.}
\bvolume{68}
\bpages{837--857}.
\bid{doi={10.1111/j.1467-9868.2006.00569.x}, issn={1369-7412}, mr={2301297}}
\bptok{imsref}%
\end{barticle}
%
\endbibitem

\bibitem[\protect\citeauthoryear{Antoniadis and
Sapatinas}{2003}]{antoniadis2003}
%
\begin{barticle}[mr]
\bauthor{\bsnm{Antoniadis},~\bfnm{Anestis}\binits{A.}} \AND
\bauthor{\bsnm{Sapatinas},~\bfnm{Theofanis}\binits{T.}}
(\byear{2003}).
\btitle{Wavelet methods for continuous-time prediction using {H}ilbert-valued
autoregressive processes}.
\bjournal{J. Multivariate Anal.}
\bvolume{87}
\bpages{133--158}.
\bid{doi={10.1016/S0047-259X(03)00028-9}, issn={0047-259X}, mr={2007265}}
\bptok{imsref}%
\end{barticle}
%
\endbibitem

\bibitem[\protect\citeauthoryear{Benko, H{\"a}rdle and
Kneip}{2009}]{benko2009}
%
\begin{barticle}[mr]
\bauthor{\bsnm{Benko},~\bfnm{Michal}\binits{M.}},
\bauthor{\bsnm{H{\"a}rdle},~\bfnm{Wolfgang}\binits{W.}} \AND
\bauthor{\bsnm{Kneip},~\bfnm{Alois}\binits{A.}}
(\byear{2009}).
\btitle{Common functional principal components}.
\bjournal{Ann. Statist.}
\bvolume{37}
\bpages{1--34}.
\bid{doi={10.1214/07-AOS516}, issn={0090-5364}, mr={2488343}}
\bptok{imsref}%
\end{barticle}
%
\endbibitem

\bibitem[\protect\citeauthoryear{Bloomfield}{2000}]{bloomfield2000}
%
\begin{bbook}[mr]
\bauthor{\bsnm{Bloomfield},~\bfnm{Peter}\binits{P.}}
(\byear{2000}).
\btitle{Fourier Analysis of Time Series: An Introduction},
\bedition{2nd} ed.
\bpublisher{Wiley}, \blocation{New York}.
\bid{doi={10.1002/0471722235}, mr={1884963}}
\bptok{imsref}%
\end{bbook}
%
\endbibitem

\bibitem[\protect\citeauthoryear{Boente, Rodriguez and
Sued}{2011}]{Boente2011}
%
\begin{bincollection}[mr]
\bauthor{\bsnm{Boente},~\bfnm{Graciela}\binits{G.}},
\bauthor{\bsnm{Rodriguez},~\bfnm{Daniela}\binits{D.}} \AND
\bauthor{\bsnm{Sued},~\bfnm{Mariela}\binits{M.}}
(\byear{2011}).
\btitle{Testing the equality of covariance operators}.
In \bbooktitle{Recent Advances in Functional Data Analysis and Related Topics}
\bpages{49--53}.
\bpublisher{Physica-Verlag/Springer}, \baddress{Heidelberg}.
\bid{doi={10.1007/978-3-7908-2736-1_8}, mr={2815560}}
\bptok{imsref}%
\end{bincollection}
%
\endbibitem

\bibitem[\protect\citeauthoryear{Bosq}{2000}]{Bosq00}
%
\begin{bbook}[mr]
\bauthor{\bsnm{Bosq},~\bfnm{D.}\binits{D.}}
(\byear{2000}).
\btitle{Linear Processes in Function Spaces: Theory and Applications}.
\bseries{Lecture Notes in Statistics}
\bvolume{149}.
\bpublisher{Springer}, \blocation{New York}.
\bid{doi={10.1007/978-1-4612-1154-9}, mr={1783138}}
\bptok{imsref}%
\end{bbook}
%
\endbibitem

\bibitem[\protect\citeauthoryear{Bosq}{2002}]{bosq2002}
%
\begin{barticle}[mr]
\bauthor{\bsnm{Bosq},~\bfnm{Denis}\binits{D.}}
(\byear{2002}).
\btitle{Estimation of mean and covariance operator of autoregressive processes
in {B}anach spaces}.
\bjournal{Stat. Inference Stoch. Process.}
\bvolume{5}
\bpages{287--306}.
\bid{doi={10.1023/A:1021279131053}, issn={1387-0874}, mr={1943835}}
\bptok{imsref}%
\end{barticle}
%
\endbibitem

\bibitem[\protect\citeauthoryear{Bosq and Blanke}{2007}]{Blanke2008}
%
\begin{bbook}[mr]
\bauthor{\bsnm{Bosq},~\bfnm{Denis}\binits{D.}} \AND
\bauthor{\bsnm{Blanke},~\bfnm{Delphine}\binits{D.}}
(\byear{2007}).
\btitle{Inference and Prediction in Large Dimensions}.
\bpublisher{Wiley}, \blocation{Chichester}.
\bid{doi={10.1002/9780470724033}, mr={2364006}}
\bptnote{check year}%
\bptok{imsref}%
\end{bbook}
%
\endbibitem

\bibitem[\protect\citeauthoryear{Brillinger}{2001}]{Brill01}
%
\begin{bbook}[mr]
\bauthor{\bsnm{Brillinger},~\bfnm{David~R.}\binits{D.~R.}}
(\byear{2001}).
\btitle{Time Series: Data Analysis and Theory}.
\bseries{Classics in Applied Mathematics}
\bvolume{36}.
\bpublisher{SIAM},
\blocation{Philadelphia, PA}.
\bid{doi={10.1137/1.9780898719246}, mr={1853554}}
\bptok{imsref}%
\end{bbook}
%
\endbibitem

\bibitem[\protect\citeauthoryear{Cardot and Sarda}{2006}]{cardot2003}
%
\begin{bincollection}[mr]
\bauthor{\bsnm{Cardot},~\bfnm{Herv{\'e}}\binits{H.}} \AND
\bauthor{\bsnm{Sarda},~\bfnm{Pascal}\binits{P.}}
(\byear{2006}).
\btitle{Linear regression models for functional data}.
In \bbooktitle{The Art of Semiparametrics}
\bpages{49--66}.
\bpublisher{Physica-Verlag/Springer}, \baddress{Heidelberg}.
\bid{doi={10.1007/3-7908-1701-5_4}, mr={2234875}}
\bptok{imsref}%
\end{bincollection}
%
\endbibitem

\bibitem[\protect\citeauthoryear{Cuevas, Febrero and
Fraiman}{2002}]{cuevas2002}
%
\begin{barticle}[mr]
\bauthor{\bsnm{Cuevas},~\bfnm{Antonio}\binits{A.}},
\bauthor{\bsnm{Febrero},~\bfnm{Manuel}\binits{M.}} \AND
\bauthor{\bsnm{Fraiman},~\bfnm{Ricardo}\binits{R.}}
(\byear{2002}).
\btitle{Linear functional regression: The case of fixed design and functional
response}.
\bjournal{Canad. J. Statist.}
\bvolume{30}
\bpages{285--300}.
\bid{doi={10.2307/3315952}, issn={0319-5724}, mr={1926066}}
\bptok{imsref}%
\end{barticle}
%
\endbibitem

\bibitem[\protect\citeauthoryear{Dauxois, Pousse and
Romain}{1982}]{dauxois1982}
%
\begin{barticle}[mr]
\bauthor{\bsnm{Dauxois},~\bfnm{J.}\binits{J.}},
\bauthor{\bsnm{Pousse},~\bfnm{A.}\binits{A.}} \AND
\bauthor{\bsnm{Romain},~\bfnm{Y.}\binits{Y.}}
(\byear{1982}).
\btitle{Asymptotic theory for the principal component analysis of a vector
random function: Some applications to statistical inference}.
\bjournal{J. Multivariate Anal.}
\bvolume{12}
\bpages{136--154}.
\bid{doi={10.1016/0047-259X(82)90088-4}, issn={0047-259X}, mr={0650934}}
\bptok{imsref}%
\end{barticle}
%
\endbibitem

\bibitem[\protect\citeauthoryear{Dehling and Sharipov}{2005}]{dehling2005}
%
\begin{barticle}[mr]
\bauthor{\bsnm{Dehling},~\bfnm{Herold}\binits{H.}} \AND
\bauthor{\bsnm{Sharipov},~\bfnm{Olimjon~Sh.}\binits{O.~S.}}
(\byear{2005}).
\btitle{Estimation of mean and covariance operator for {B}anach space valued
autoregressive processes with dependent innovations}.
\bjournal{Stat. Inference Stoch. Process.}
\bvolume{8}
\bpages{137--149}.
\bid{doi={10.1007/s11203-003-0382-8}, issn={1387-0874}, mr={2121674}}
\bptok{imsref}%
\end{barticle}
%
\endbibitem

\bibitem[\protect\citeauthoryear{Edwards}{1967}]{Edwards1967}
%
\begin{bbook}[auto:STB|2013/03/04|13:35:07]
\bauthor{\bsnm{Edwards},~\bfnm{R.}\binits{R.}}
(\byear{1967}).
\btitle{Fourier Series: A Modern Introduction}.
\bpublisher{Holt, Rinehart \& Winston}, \blocation{New York}.
\bptok{imsref}%
\end{bbook}
%
\endbibitem

\bibitem[\protect\citeauthoryear{Ferraty and Vieu}{2004}]{ferraty2004}
%
\begin{barticle}[mr]
\bauthor{\bsnm{Ferraty},~\bfnm{F.}\binits{F.}} \AND
\bauthor{\bsnm{Vieu},~\bfnm{P.}\binits{P.}}
(\byear{2004}).
\btitle{Nonparametric models for functional data, with application in
regression, time-series prediction and curve discrimination}.
\bjournal{J. Nonparametr. Stat.}
\bvolume{16}
\bpages{111--125}.
\bid{doi={10.1080/10485250310001622686}, issn={1048-5252}, mr={2053065}}
\bptok{imsref}%
\end{barticle}
%
\endbibitem

\bibitem[\protect\citeauthoryear{Ferraty and Vieu}{2006}]{ferraty2006}
%
\begin{bbook}[mr]
\bauthor{\bsnm{Ferraty},~\bfnm{Fr{\'e}d{\'e}ric}\binits{F.}} \AND
\bauthor{\bsnm{Vieu},~\bfnm{Philippe}\binits{P.}}
(\byear{2006}).
\btitle{Nonparametric Functional Data Analysis: Theory and Practice}.
\bpublisher{Springer}, \blocation{New York}.
\bid{mr={2229687}}
\bptok{imsref}%
\end{bbook}
%
\endbibitem

\bibitem[\protect\citeauthoryear{Ferraty et~al.}{2011a}]{ferraty-goia2011}
%
\begin{bincollection}[mr]
\bauthor{\bsnm{Ferraty},~\bfnm{Fr{\'e}d{\'e}ric}\binits{F.}},
\bauthor{\bsnm{Goia},~\bfnm{Aldo}\binits{A.}},
\bauthor{\bsnm{Salinelli},~\bfnm{Enersto}\binits{E.}} \AND
\bauthor{\bsnm{Vieu},~\bfnm{Philippe}\binits{P.}}
(\byear{2011}a).
\btitle{Recent advances on functional additive regression}.
In \bbooktitle{Recent Advances in Functional Data Analysis and Related Topics}
\bpages{97--102}.
\bpublisher{Physica-Verlag/Springer}, \baddress{Heidelberg}.
\bid{doi={10.1007/978-3-7908-2736-1_15}, mr={2815567}}
\bptok{imsref}%
\end{bincollection}
%
\endbibitem

\bibitem[\protect\citeauthoryear{Ferraty et~al.}{2011b}]{ferraty-tadj2011}
%
\begin{barticle}[mr]
\bauthor{\bsnm{Ferraty},~\bfnm{Fr{\'e}d{\'e}ric}\binits{F.}},
\bauthor{\bsnm{Laksaci},~\bfnm{Ali}\binits{A.}},
\bauthor{\bsnm{Tadj},~\bfnm{Amel}\binits{A.}} \AND
\bauthor{\bsnm{Vieu},~\bfnm{Philippe}\binits{P.}}
(\byear{2011}b).
\btitle{Kernel regression with functional response}.
\bjournal{Electron. J. Stat.}
\bvolume{5}
\bpages{159--171}.
\bid{doi={10.1214/11-EJS600}, issn={1935-7524}, mr={2786486}}
\bptok{imsref}%
\end{barticle}
%
\endbibitem

\bibitem[\protect\citeauthoryear{Fremdt et~al.}{2013}]{Fremdt2011}
%
\begin{barticle}[auto:STB|2013/03/04|13:35:07]
\bauthor{\bsnm{Fremdt},~\bfnm{S.}\binits{S.}},
\bauthor{\bsnm{Steinebach},~\bfnm{J.}\binits{J.}},
\bauthor{\bsnm{Horv{\'a}th},~\bfnm{L.}\binits{L.}}
\AND
\bauthor{\bsnm{Kokoszka},~\bfnm{P.}\binits{P.}}
(\byear{2013}).
\btitle{Testing the equality of covariance operators in functional
samples}.
\bjournal{Scand. J. Stat.}
\bvolume{40}
\bpages{138--152}.
\bptok{imsref}%
\end{barticle}
%
\endbibitem

\bibitem[\protect\citeauthoryear{Gabrys, Horv{\'a}th and
Kokoszka}{2010}]{Gabrys2010}
%
\begin{barticle}[mr]
\bauthor{\bsnm{Gabrys},~\bfnm{Robertas}\binits{R.}},
\bauthor{\bsnm{Horv{\'a}th},~\bfnm{Lajos}\binits{L.}} \AND
\bauthor{\bsnm{Kokoszka},~\bfnm{Piotr}\binits{P.}}
(\byear{2010}).
\btitle{Tests for error correlation in the functional linear model}.
\bjournal{J. Amer. Statist. Assoc.}
\bvolume{105}
\bpages{1113--1125}.
\bid{doi={10.1198/jasa.2010.tm09794}, issn={0162-1459}, mr={2752607}}
\bptok{imsref}%
\end{barticle}
%
\endbibitem

\bibitem[\protect\citeauthoryear{Gabrys and Kokoszka}{2007}]{Gabrys2007}
%
\begin{barticle}[mr]
\bauthor{\bsnm{Gabrys},~\bfnm{Robertas}\binits{R.}} \AND
\bauthor{\bsnm{Kokoszka},~\bfnm{Piotr}\binits{P.}}
(\byear{2007}).
\btitle{Portmanteau test of independence for functional observations}.
\bjournal{J. Amer. Statist. Assoc.}
\bvolume{102}
\bpages{1338--1348}.
\bid{doi={10.1198/016214507000001111}, issn={0162-1459}, mr={2412554}}
\bptok{imsref}%
\end{barticle}
%
\endbibitem

\bibitem[\protect\citeauthoryear{Grenander}{1981}]{grenander1981}
%
\begin{bbook}[mr]
\bauthor{\bsnm{Grenander},~\bfnm{Ulf}\binits{U.}}
(\byear{1981}).
\btitle{Abstract Inference}.
\bpublisher{Wiley}, \blocation{New York}.
\bid{mr={0599175}}
\bptok{imsref}%
\end{bbook}
%
\endbibitem

\bibitem[\protect\citeauthoryear{Grenander and
Rosenblatt}{1957}]{Grenander1966}
%
\begin{bbook}[mr]
\bauthor{\bsnm{Grenander},~\bfnm{Ulf}\binits{U.}} \AND
\bauthor{\bsnm{Rosenblatt},~\bfnm{Murray}\binits{M.}}
(\byear{1957}).
\btitle{Statistical Analysis of Stationary Time Series}.
\bpublisher{Wiley}, \blocation{New York}.
\bid{mr={0084975}}
\bptnote{check year}%
\bptok{imsref}%
\end{bbook}
%
\endbibitem

\bibitem[\protect\citeauthoryear{Hall and Hosseini-Nasab}{2006}]{hall2006}
%
\begin{barticle}[mr]
\bauthor{\bsnm{Hall},~\bfnm{Peter}\binits{P.}} \AND
\bauthor{\bsnm{Hosseini-Nasab},~\bfnm{Mohammad}\binits{M.}}
(\byear{2006}).
\btitle{On properties of functional principal components analysis}.
\bjournal{J. R. Stat. Soc. Ser. B Stat. Methodol.}
\bvolume{68}
\bpages{109--126}.
\bid{doi={10.1111/j.1467-9868.2005.00535.x}, issn={1369-7412}, mr={2212577}}
\bptok{imsref}%
\end{barticle}
%
\endbibitem

\bibitem[\protect\citeauthoryear{Hall and Vial}{2006}]{hall2006b}
%
\begin{barticle}[mr]
\bauthor{\bsnm{Hall},~\bfnm{Peter}\binits{P.}} \AND
\bauthor{\bsnm{Vial},~\bfnm{C{\'e}line}\binits{C.}}
(\byear{2006}).
\btitle{Assessing the finite dimensionality of functional data}.
\bjournal{J.~R.~Stat. Soc. Ser. B Stat. Methodol.}
\bvolume{68}
\bpages{689--705}.
\bid{doi={10.1111/j.1467-9868.2006.00562.x}, issn={1369-7412}, mr={2301015}}
\bptok{imsref}%
\end{barticle}
%
\endbibitem

\bibitem[\protect\citeauthoryear{Hannan}{1970}]{Hannan1970}
%
\begin{bbook}[mr]
\bauthor{\bsnm{Hannan},~\bfnm{E.~J.}\binits{E.~J.}}
(\byear{1970}).
\btitle{Multiple Time Series}.
\bpublisher{Wiley}, \blocation{New York}.
\bid{mr={0279952}}
\bptok{imsref}%
\end{bbook}
%
\endbibitem

\bibitem[\protect\citeauthoryear{H{\"o}rmann and Kokoszka}{2010}]{Hormann10}
%
\begin{barticle}[mr]
\bauthor{\bsnm{H{\"o}rmann},~\bfnm{Siegfried}\binits{S.}} \AND
\bauthor{\bsnm{Kokoszka},~\bfnm{Piotr}\binits{P.}}
(\byear{2010}).
\btitle{Weakly dependent functional data}.
\bjournal{Ann. Statist.}
\bvolume{38}
\bpages{1845--1884}.
\bid{doi={10.1214/09-AOS768}, issn={0090-5364}, mr={2662361}}
\bptok{imsref}%
\end{barticle}
%
\endbibitem

\bibitem[\protect\citeauthoryear{Horv{\'a}th, Hu{\v{s}}kov{\'a} and
Kokoszka}{2010}]{Horvath10a}
%
\begin{barticle}[mr]
\bauthor{\bsnm{Horv{\'a}th},~\bfnm{Lajos}\binits{L.}},
\bauthor{\bsnm{Hu{\v{s}}kov{\'a}},~\bfnm{Marie}\binits{M.}} \AND
\bauthor{\bsnm{Kokoszka},~\bfnm{Piotr}\binits{P.}}
(\byear{2010}).
\btitle{Testing the stability of the functional autoregressive process}.
\bjournal{J. Multivariate Anal.}
\bvolume{101}
\bpages{352--367}.
\bid{doi={10.1016/j.jmva.2008.12.008}, issn={0047-259X}, mr={2564345}}
\bptok{imsref}%
\end{barticle}
%
\endbibitem

\bibitem[\protect\citeauthoryear{Horv{\'a}th and Kokoszka}{2012}]{Horvath2012}
%
\begin{bbook}[mr]
\bauthor{\bsnm{Horv{\'a}th},~\bfnm{Lajos}\binits{L.}} \AND
\bauthor{\bsnm{Kokoszka},~\bfnm{Piotr}\binits{P.}}
(\byear{2012}).
\btitle{Inference for Functional Data with Applications}.
\bpublisher{Springer}, \blocation{New York}.
\bid{doi={10.1007/978-1-4614-3655-3}, mr={2920735}}
\bptok{imsref}%
\end{bbook}
%
\endbibitem

\bibitem[\protect\citeauthoryear{Horv{\'a}th, Kokoszka and
Reeder}{2013}]{Horvath11}
%
\begin{barticle}[auto:STB|2013/03/04|13:35:07]
\bauthor{\bsnm{Horv{\'a}th},~\bfnm{L.}\binits{L.}},
\bauthor{\bsnm{Kokoszka},~\bfnm{P.}\binits{P.}} \AND
\bauthor{\bsnm{Reeder},~\bfnm{R.}\binits{R.}}
(\byear{2013}).
\btitle{Estimation of the mean of functional time
series and a
two-sample problem.}
\bjournal{J. R. Stat. Soc. Ser. B Stat. Methodol.}
\bvolume{75}
\bpages{103--122}.
\bptok{imsref}%
\end{barticle}
%
\endbibitem

\bibitem[\protect\citeauthoryear{Hunter and Nachtergaele}{2001}]{Hunter2005}
%
\begin{bbook}[mr]
\bauthor{\bsnm{Hunter},~\bfnm{John~K.}\binits{J.~K.}} \AND
\bauthor{\bsnm{Nachtergaele},~\bfnm{Bruno}\binits{B.}}
(\byear{2001}).
\btitle{Applied Analysis}.
\bpublisher{World Scientific}, \blocation{River Edge, NJ}.
\bid{mr={1829589}}
\bptnote{check year}%
\bptok{imsref}%
\end{bbook}
%
\endbibitem

\bibitem[\protect\citeauthoryear{Kadison and Ringrose}{1997}]{Kadison1997}
%
\begin{bbook}[auto:STB|2013/03/04|13:35:07]
\bauthor{\bsnm{Kadison},~\bfnm{R.~V.}\binits{R.~V.}} \AND
\bauthor{\bsnm{Ringrose},~\bfnm{J.~R.}\binits{J.~R.}}
(\byear{1997}).
\btitle{Fundamentals of the Theory of Operator Algebras}.
\bseries{Graduate Studies in Mathematics}
\bvolume{15}.
\bpublisher{Amer. Math. Soc.}, \blocation{Providence, RI}.
\bptok{imsref}%
\end{bbook}
%
\endbibitem

\bibitem[\protect\citeauthoryear{Karhunen}{1947}]{karhunen1947}
%
\begin{barticle}[mr]
\bauthor{\bsnm{Karhunen},~\bfnm{Kari}\binits{K.}}
(\byear{1947}).
\btitle{\"{U}ber lineare {M}ethoden in der {W}ahrscheinlichkeitsrechnung}.
\bjournal{Ann. Acad. Sci. Fennicae. Ser. A I Math.-Phys.}
\bvolume{1947}
\bpages{79}.
\bid{mr={0023013}}
\bptok{imsref}%
\end{barticle}
%
\endbibitem

\bibitem[\protect\citeauthoryear{Kolmogorov}{1978}]{kolmogorov1978}
%
\begin{bbook}[auto:STB|2013/03/04|13:35:07]
\bauthor{\bsnm{Kolmogorov},~\bfnm{A.}\binits{A.}}
(\byear{1978}).
\btitle{Stationary Sequences in Hilbert Space}.
\bpublisher{National Translations Center [John Crerar Library]}, \blocation{Chicago}.
\bptok{imsref}%
\end{bbook}
%
\endbibitem

\bibitem[\protect\citeauthoryear{Kraus and Panaretos}{2012}]{kraus2012}
%
\begin{barticle}[auto:STB|2013/03/04|13:35:07]
\bauthor{\bsnm{Kraus},~\bfnm{D.}\binits{D.}} \AND
\bauthor{\bsnm{Panaretos},~\bfnm{V.~M.}\binits{V.~M.}}
(\byear{2012}).
\btitle{Disperson operators and resistant second-order functional data
analysis}.
\bjournal{Biometrika}
\bvolume{99}
\bpages{813--832}.
\bptok{imsref}%
\end{barticle}
%
\endbibitem

\bibitem[\protect\citeauthoryear{Laib and Louani}{2010}]{louani2010}
%
\begin{barticle}[mr]
\bauthor{\bsnm{Laib},~\bfnm{Na{\^a}mane}\binits{N.}} \AND
\bauthor{\bsnm{Louani},~\bfnm{Djamal}\binits{D.}}
(\byear{2010}).
\btitle{Nonparametric kernel regression estimation for functional stationary
ergodic data: Asymptotic properties}.
\bjournal{J. Multivariate Anal.}
\bvolume{101}
\bpages{2266--2281}.
\bid{doi={10.1016/j.jmva.2010.05.010}, issn={0047-259X}, mr={2719861}}
\bptok{imsref}%
\end{barticle}
%
\endbibitem

\bibitem[\protect\citeauthoryear{Ledoux and Talagrand}{1991}]{Ledo91}
%
\begin{bbook}[mr]
\bauthor{\bsnm{Ledoux},~\bfnm{Michel}\binits{M.}} \AND
\bauthor{\bsnm{Talagrand},~\bfnm{Michel}\binits{M.}}
(\byear{1991}).
\btitle{Probability in {B}anach Spaces: Isoperimetry and Processes}.
\bseries{Ergebnisse der Mathematik und Ihrer Grenzgebiete (3) [Results in
Mathematics and Related Areas (3)]}
\bvolume{23}.
\bpublisher{Springer}, \blocation{Berlin}.
\bid{mr={1102015}}
\bptok{imsref}%
\end{bbook}
%
\endbibitem

\bibitem[\protect\citeauthoryear{L{\'e}vy}{1948}]{loeve1948}
%
\begin{bbook}[mr]
\bauthor{\bsnm{L{\'e}vy},~\bfnm{Paul}\binits{P.}}
(\byear{1948}).
\btitle{Processus {s}tochastiques et {m}ouvement {B}rownien. {S}uivi
d'une note
de {M}.~{L}o\`eve}.
\bpublisher{Gauthier-Villars}, \blocation{Paris}.
\bid{mr={0029120}}
\bptok{imsref}%
\end{bbook}
%
\endbibitem

\bibitem[\protect\citeauthoryear{Liu and Wu}{2010}]{liu-wu2010}
%
\begin{barticle}[mr]
\bauthor{\bsnm{Liu},~\bfnm{Weidong}\binits{W.}} \AND
\bauthor{\bsnm{Wu},~\bfnm{Wei~Biao}\binits{W.~B.}}
(\byear{2010}).
\btitle{Asymptotics of spectral density estimates}.
\bjournal{Econometric Theory}
\bvolume{26}
\bpages{1218--1245}.
\bid{doi={10.1017/S026646660999051X}, issn={0266-4666}, mr={2660298}}
\bptok{imsref}%
\end{barticle}
%
\endbibitem

\bibitem[\protect\citeauthoryear{Locantore et~al.}{1999}]{locantore1999}
%
\begin{barticle}[mr]
\bauthor{\bsnm{Locantore},~\bfnm{N.}\binits{N.}},
\bauthor{\bsnm{Marron},~\bfnm{J.~S.}\binits{J.~S.}},
\bauthor{\bsnm{Simpson},~\bfnm{D.~G.}\binits{D.~G.}},
\bauthor{\bsnm{Tripoli},~\bfnm{N.}\binits{N.}},
\bauthor{\bsnm{Zhang},~\bfnm{J.~T.}\binits{J.~T.}} \AND
\bauthor{\bsnm{Cohen},~\bfnm{K.~L.}\binits{K.~L.}}
(\byear{1999}).
\btitle{Robust principal component analysis for functional data}.
\bjournal{TEST}
\bvolume{8}
\bpages{1--73}.
\bid{doi={10.1007/BF02595862}, issn={1133-0686}, mr={1707596}}
\bptnote{check related}%
\bptok{imsref}%
\end{barticle}
%
\endbibitem

\bibitem[\protect\citeauthoryear{Mas}{2000}]{mas2000}
%
\begin{bmisc}[auto:STB|2013/03/04|13:35:07]
\bauthor{\bsnm{Mas},~\bfnm{A.}\binits{A.}}
(\byear{2000}).
\bhowpublished{Estimation d'op\'erateurs de corr\'elation de processus
lin\'eaires fonctionnels: lois limites, d\'eviations mod\'er\'ees.
Ph.D. thesis, Universit\'e Paris VI}.
\bptok{imsref}%
\end{bmisc}
%
\endbibitem

\bibitem[\protect\citeauthoryear{Panaretos, Kraus and
Maddocks}{2010}]{panaretos2010}
%
\begin{barticle}[mr]
\bauthor{\bsnm{Panaretos},~\bfnm{Victor~M.}\binits{V.~M.}},
\bauthor{\bsnm{Kraus},~\bfnm{David}\binits{D.}} \AND
\bauthor{\bsnm{Maddocks},~\bfnm{John~H.}\binits{J.~H.}}
(\byear{2010}).
\btitle{Second-order comparison of {G}aussian random functions and the geometry
of {DNA} minicircles}.
\bjournal{J. Amer. Statist. Assoc.}
\bvolume{105}
\bpages{670--682}.
\bid{doi={10.1198/jasa.2010.tm09239}, issn={0162-1459}, mr={2724851}}
\bptok{imsref}%
\end{barticle}
%
\endbibitem

\bibitem[\protect\citeauthoryear{Panaretos and Tavakoli}{2013}]{stoc}
%
\begin{bmisc}[auto:STB|2013/03/04|13:35:07]
\bauthor{\bsnm{Panaretos},~\bfnm{V.~M.}\binits{V.~M.}} \AND
\bauthor{\bsnm{Tavakoli},~\bfnm{S.}\binits{S.}}
(\byear{2013}).
\bhowpublished{Cram\'er--Karhunen--Lo\`eve representation and harmonic principal
component analysis of functional time series. \textit{Stochastic Process. Appl.} To appear.
DOI:\doiurl{10.1016/j.spa.2013.03.015}, available at
\url{http://www.sciencedirect.com/science/article/pii/S0304414913000793}.}
\bptok{imsref}%
\end{bmisc}
%
\endbibitem

\bibitem[\protect\citeauthoryear{Panaretos and Tavakoli}{2013}]{supp}
%
\begin{bmisc}[auto]
\bauthor{\bsnm{Panaretos},~\bfnm{V.~M.}\binits{V.~M.}} \AND
\bauthor{\bsnm{Tavakoli},~\bfnm{S.}\binits{S.}}
(\byear{2013}).
\bhowpublished{Supplement to ``Fourier analysis of stationary time series in function
space.'' DOI:\doiurl{10.1214/13-AOS1086SUPP}.}
\bptok{imsref}%
\end{bmisc}
%
\endbibitem


\bibitem[\protect\citeauthoryear{Peligrad and Wu}{2010}]{peligrad-wu2010}
%
\begin{barticle}[mr]
\bauthor{\bsnm{Peligrad},~\bfnm{Magda}\binits{M.}} \AND
\bauthor{\bsnm{Wu},~\bfnm{Wei~Biao}\binits{W.~B.}}
(\byear{2010}).
\btitle{Central limit theorem for {F}ourier transforms of stationary
processes}.
\bjournal{Ann. Probab.}
\bvolume{38}
\bpages{2009--2022}.
\bid{doi={10.1214/10-AOP530}, issn={0091-1798}, mr={2722793}}
\bptok{imsref}%
\end{barticle}
%
\endbibitem

\bibitem[\protect\citeauthoryear{Pollard}{1984}]{Pollard1984}
%
\begin{bbook}[mr]
\bauthor{\bsnm{Pollard},~\bfnm{David}\binits{D.}}
(\byear{1984}).
\btitle{Convergence of Stochastic Processes}.
\bpublisher{Springer}, \blocation{New York}.
\bid{doi={10.1007/978-1-4612-5254-2}, mr={0762984}}
\bptok{imsref}%
\end{bbook}
%
\endbibitem

\bibitem[\protect\citeauthoryear{Priestley}{2001}]{priestley2001}
%
\begin{bbook}[auto:STB|2013/03/04|13:35:07]
\bauthor{\bsnm{Priestley},~\bfnm{M.~B.}\binits{M.~B.}}
(\byear{2001}).
\btitle{Spectral Analysis and Time Series, Vol. I and II}.
\bpublisher{Academic Press}, \blocation{San Diego}.
\bptok{imsref}%
\end{bbook}
%
\endbibitem

\bibitem[\protect\citeauthoryear{Ramsay and Silverman}{2005}]{ramsay2005}
%
\begin{bbook}[mr]
\bauthor{\bsnm{Ramsay},~\bfnm{J.~O.}\binits{J.~O.}} \AND
\bauthor{\bsnm{Silverman},~\bfnm{B.~W.}\binits{B.~W.}}
(\byear{2005}).
\btitle{Functional Data Analysis},
\bedition{2nd} ed.
\bpublisher{Springer}, \blocation{New York}.
\bid{mr={2168993}}
\bptok{imsref}%
\end{bbook}
%
\endbibitem

\bibitem[\protect\citeauthoryear{Rice and Silverman}{1991}]{rice1991}
%
\begin{barticle}[mr]
\bauthor{\bsnm{Rice},~\bfnm{John~A.}\binits{J.~A.}} \AND
\bauthor{\bsnm{Silverman},~\bfnm{B.~W.}\binits{B.~W.}}
(\byear{1991}).
\btitle{Estimating the mean and covariance structure nonparametrically
when the
data are curves}.
\bjournal{J. Roy. Statist. Soc. Ser. B}
\bvolume{53}
\bpages{233--243}.
\bid{issn={0035-9246}, mr={1094283}}
\bptok{imsref}%
\end{barticle}
%
\endbibitem

\bibitem[\protect\citeauthoryear{Rosenblatt}{1984}]{Rosenblatt1984}
%
\begin{barticle}[mr]
\bauthor{\bsnm{Rosenblatt},~\bfnm{M.}\binits{M.}}
(\byear{1984}).
\btitle{Asymptotic normality, strong mixing and spectral density estimates}.
\bjournal{Ann. Probab.}
\bvolume{12}
\bpages{1167--1180}.
\bid{issn={0091-1798}, mr={0757774}}
\bptok{imsref}%
\end{barticle}
%
\endbibitem

\bibitem[\protect\citeauthoryear{Rosenblatt}{1985}]{Rosenblatt1985}
%
\begin{bbook}[mr]
\bauthor{\bsnm{Rosenblatt},~\bfnm{Murray}\binits{M.}}
(\byear{1985}).
\btitle{Stationary Sequences and Random Fields}.
\bpublisher{Birkh\"auser}, \blocation{Boston, MA}.
\bid{doi={10.1007/978-1-4612-5156-9}, mr={0885090}}
\bptok{imsref}%
\end{bbook}
%
\endbibitem

\bibitem[\protect\citeauthoryear{Sen and Kl{\"u}ppelberg}{2010}]{Sen10}
%
\begin{bmisc}[auto:STB|2013/03/04|13:35:07]
\bauthor{\bsnm{Sen},~\bfnm{R.}\binits{R.}} \AND
\bauthor{\bsnm{Kl{\"u}ppelberg},~\bfnm{C.}\binits{C.}}
(\byear{2010}).
\bhowpublished{Time series of functional data. Unpublished manuscript.
Available at
\texttt{\href{http://citeseerx.ist.psu.edu/viewdoc/summary?doi=10.1.1.185.2739}{http://citeseerx.ist.psu.edu/viewdoc/summary?doi=}
\href{http://citeseerx.ist.psu.edu/viewdoc/summary?doi=10.1.1.185.2739}{10.1.1.185.2739}}.}
\bptok{imsref}%
\end{bmisc}
%
\endbibitem

\bibitem[\protect\citeauthoryear{Shao and Wu}{2007}]{shao-wu2007}
%
\begin{barticle}[mr]
\bauthor{\bsnm{Shao},~\bfnm{Xiaofeng}\binits{X.}} \AND
\bauthor{\bsnm{Wu},~\bfnm{Wei~Biao}\binits{W.~B.}}
(\byear{2007}).
\btitle{Asymptotic spectral theory for nonlinear time series}.
\bjournal{Ann. Statist.}
\bvolume{35}
\bpages{1773--1801}.
\bid{doi={10.1214/009053606000001479}, issn={0090-5364}, mr={2351105}}
\bptok{imsref}%
\end{barticle}
%
\endbibitem

\bibitem[\protect\citeauthoryear{Wand and Jones}{1995}]{Wand1995}
%
\begin{bbook}[mr]
\bauthor{\bsnm{Wand},~\bfnm{M.~P.}\binits{M.~P.}} \AND
\bauthor{\bsnm{Jones},~\bfnm{M.~C.}\binits{M.~C.}}
(\byear{1995}).
\btitle{Kernel Smoothing}.
\bseries{Monographs on Statistics and Applied Probability}
\bvolume{60}.
\bpublisher{Chapman \& Hall}, \blocation{London}.
\bid{mr={1319818}}
\bptok{imsref}%
\end{bbook}
%
\endbibitem

\bibitem[\protect\citeauthoryear{Weidmann}{1980}]{Weid80}
%
\begin{bbook}[mr]
\bauthor{\bsnm{Weidmann},~\bfnm{Joachim}\binits{J.}}
(\byear{1980}).
\btitle{Linear Operators in {H}ilbert Spaces}.
\bseries{Graduate Texts in Mathematics}
\bvolume{68}.
\bpublisher{Springer}, \blocation{New York}.
\bid{mr={0566954}}
\bptok{imsref}%
\end{bbook}
%
\endbibitem

\bibitem[\protect\citeauthoryear{Wheeden and Zygmund}{1977}]{Wheeden1977}
%
\begin{bbook}[mr]
\bauthor{\bsnm{Wheeden},~\bfnm{Richard~L.}\binits{R.~L.}} \AND
\bauthor{\bsnm{Zygmund},~\bfnm{Antoni}\binits{A.}}
(\byear{1977}).
\btitle{Measure and Integral: An Introduction to Real Analysis}.
\bseries{Pure and Applied Mathematics}
\bvolume{43}.
\bpublisher{Dekker}, \blocation{New York}.
\bid{mr={0492146}}
\bptok{imsref}%
\end{bbook}
%
\endbibitem

\bibitem[\protect\citeauthoryear{Yao, M{\"u}ller and Wang}{2005}]{yao2005}
%
\begin{barticle}[mr]
\bauthor{\bsnm{Yao},~\bfnm{Fang}\binits{F.}},
\bauthor{\bsnm{M{\"u}ller},~\bfnm{Hans-Georg}\binits{H.-G.}} \AND
\bauthor{\bsnm{Wang},~\bfnm{Jane-Ling}\binits{J.-L.}}
(\byear{2005}).
\btitle{Functional linear regression analysis for longitudinal data}.
\bjournal{Ann. Statist.}
\bvolume{33}
\bpages{2873--2903}.
\bid{doi={10.1214/009053605000000660}, issn={0090-5364}, mr={2253106}}
\bptok{imsref}%
\end{barticle}
%
\endbibitem

\end{thebibliography}

%

\printaddresses

\end{document}